\documentclass{rspublic}
\pdfoutput=1
\usepackage[numbers]{natbib}
\usepackage[font=footnotesize]{caption}
\usepackage{color,graphicx,amsmath,gensymb}
\usepackage[normalem]{ulem}
\usepackage{mathptmx}
\renewcommand{\d}{\mathop{}\!\mathrm{d}}
\renewcommand{\epsilon}{\varepsilon}
\newcommand{\var}{\operatorname{var}}
\newcommand{\jsmod}[1]{{#1}}
\newcommand{\jsrem}[1]{}
\newcommand{\lenton}{[Lenton \emph{et al.} this volume]\ }
\newcommand{\beaulieu}{[Beaulieu \emph{et al.} this volume]\ }
\newcommand{\ashwin}{[Ashwin \emph{et al.} this volume]\ }

\title{Nonlinear softening as a predictive precursor to climate tipping}

\author[J. Sieber and J. M.T. Thompson]{Jan Sieber$^1$ and J. Michael
  T. Thompson$^2$}

 \affiliation{
   $^1$Department of Mathematics, University of Portsmouth, UK\\
 $^2$School of Engineering, University of
      Aberdeen, UK,  and
      Department of Applied Mathematics and Theoretical Physics,
      University of Cambridge, UK}
\begin{document}
\maketitle
\begin{abstract}{time series analysis, bifurcation prediction, climate tipping}
  Approaching a dangerous bifurcation, from which a dynamical system
  such as the Earth's climate will jump (tip) to a different state,
  the current stable state lies within a shrinking basin of
  attraction. Persistence of the state becomes increasingly precarious
  in the presence of noisy disturbances.  \jsrem{We consider an
    underlying potential, as defined theoretically for a saddle-node
    fold and (via averaging) for a Hopf bifurcation. Close to a stable
    state, this potential has a parabolic form; but approaching a jump
    it becomes increasingly dominated by softening nonlinearities. If
    we have already detected a decrease in the linear decay rate,
    nonlinear information allows us to estimate the propensity for
    early tipping due to noise.} \jsmod{ We argue that one needs to
    extract information about the nonlinear features (a ``softening'')
    of the underlying potential from the time series to judge the
    probability and timing of tipping. This analysis is the logical
    next step if one has detected a decrease of the linear decay
    rate.} If there is no discernable trend in the linear analysis,
  nonlinear softening is even more important in showing the proximity
  to tipping.

  After extensive normal form calibration studies, we \jsmod{ check
    two geological time series from paleo-climate tipping events for
    softening of the underlying well}. For the ending of the last ice
  age, where we find no convincing linear precursor, we identify a
  statistically significant nonlinear softening towards increasing
  temperature. \jsrem{The analysis has thus successfully detected a
    warning of the imminent tipping event.}
\end{abstract}
\section{Introduction}
\label{sec:intro}
The on-going work by the United Nations, following the major report of
the \citet{IPCC2007} in 2007 and the Climate Change Conferences in
Copenhagen and Cancun (Mexico), centres on the prediction of future
climate change, a key feature of which would be any suggestion of a
sudden and (possibly) irreversible abrupt change called a tipping
point (\citet{Lenton2008,Scheffer2009a}). Many tipping points, such as
the switching on and off of ice-ages, are well documented in
paleo-climate studies over millions of years of the Earth's history.

Predicting such tippings in advance using time-series data derived
from preceding behaviour is now seen as a major challenge, impinging,
for example on the possible use of geo-engineering
(\citet{Launder2010}). Techniques introduced by \citet{Held2004} and
\citet{Livina2007} draw on the assumption that the tipping events are
governed by a bifurcation in an underlying nonlinear dissipative
dynamical system. Specifically, researchers search for a slowing down
of intrinsic transient responses within the data, which is predicted
to occur before most bifurcational instabilities. This is done by
determining a so-called propagator, estimated from the correlation
between successive elements of a window sliding along the time series
\jsmod{(this estimate is called $\kappa_\mathrm{ACF}$ in this
  paper)}. This propagator, such as the AR(1) mapping coefficient, is
a measure of the linear decay rate and should increase to unity at a
tipping instability (corresponding to a decrease of the linear decay
rate to zero).

Prediction techniques can be tested on climatic computer models, but
more challenging is to try to predict real ancient climate tippings,
using their preceding geological data provided by ice cores,
sediments, etc. Using this preceding data alone, the aim would be to
see to what extent the actual tipping could have been predicted in
advance. 

One such study by \citet{Livina2007} looks at the end of the most
recent ice age and the associated Younger Dryas event, about $11,500$
years ago, when the Arctic warmed by $7\celsius$ in $50$ years. This
pioneering study used a time series derived from Greenland ice-core
paleo-temperature data. A second such study (one of eight made by
\citet{Dakos2008}), using data from tropical Pacific sediments, gives
a good prediction for the end of `greenhouse' Earth about $34$ million
years ago when the climate tipped from a tropical state into an
icehouse state.  \lenton gives a complete overview of
current techniques for extracting early-warning signals from time
series and compares them on realistic models and a range of
paleo-climate time series.

A \jsmod{remark} that is made during the analysis of \lenton is
that the absolute values of the extracted quantities have no direct
bearing on the tipping probability over time unless one knows also the
size of the \emph{basin of attraction}. This basin is determined in
the simplest cases by the dominant nonlinear term in the underlying
equations of motion, which is the subject of this paper.

\section{Concepts from Nonlinear Dynamics}
\label{sec:concepts}
The core techniques for analysis of time series of climate data aim to
extract the linear decay rate toward (a possibly drifting)
noise-perturbed equilibrium \citep{Held2004,Livina2007}. These
techniques can be directed at the Earth's climate as a whole, or at
the relevant climate sub-systems described by \citet{Lenton2008} as
\emph{tipping elements}. Our aim here is to augment this linear information
with information about nonlinear features of the underlying dynamics,
by examining to what extent we can identify nonlinear softening and
include it into the list of early-warning signs for climate tipping.

\subsection{Shrinking basins around dangerous bifurcations}
\label{sec:basins}
As we approach a dangerous bifurcation, from which a nonlinear
dissipative dynamical system will experience a finite jump to a remote
alternative state, the current attractor is located within a shrinking
basin of attraction. Correspondingly, the maintenance of the state
will become increasingly precarious in the presence of noisy
disturbances.

Thinking in terms of an underlying potential energy function (the
existence of which is theoretically well-defined for a saddle-node
fold and some other bifurcations) the parabolic shape of its graph
will become increasingly perturbed by softening nonlinear features.

% We show in this paper how these features can be detected and displayed
% to give a clear warning of an imminent tipping event. Specifically, if
% we have already detected a decrease in the linear decay rate, the
% nonlinear information allows us to estimate the propensity for early
% escape due to noise. If there is no discernable trend in the linear
% decay rate the nonlinear information is even more important in showing
% the (seemingly missed) proximity to tipping. \citet{Livina2010} used
% nonlinear information to detect multiple-well potentials in
% time-series. However, the study by \citet{Livina2010} was about
% \emph{identification} instead of \emph{prediction} as the time series
% windows included the transitions. See also \citet{Beaulieu2011} for
% methods that detect change points in time series.

\begin{figure}[th]
  \centering
  \includegraphics[scale=0.3]{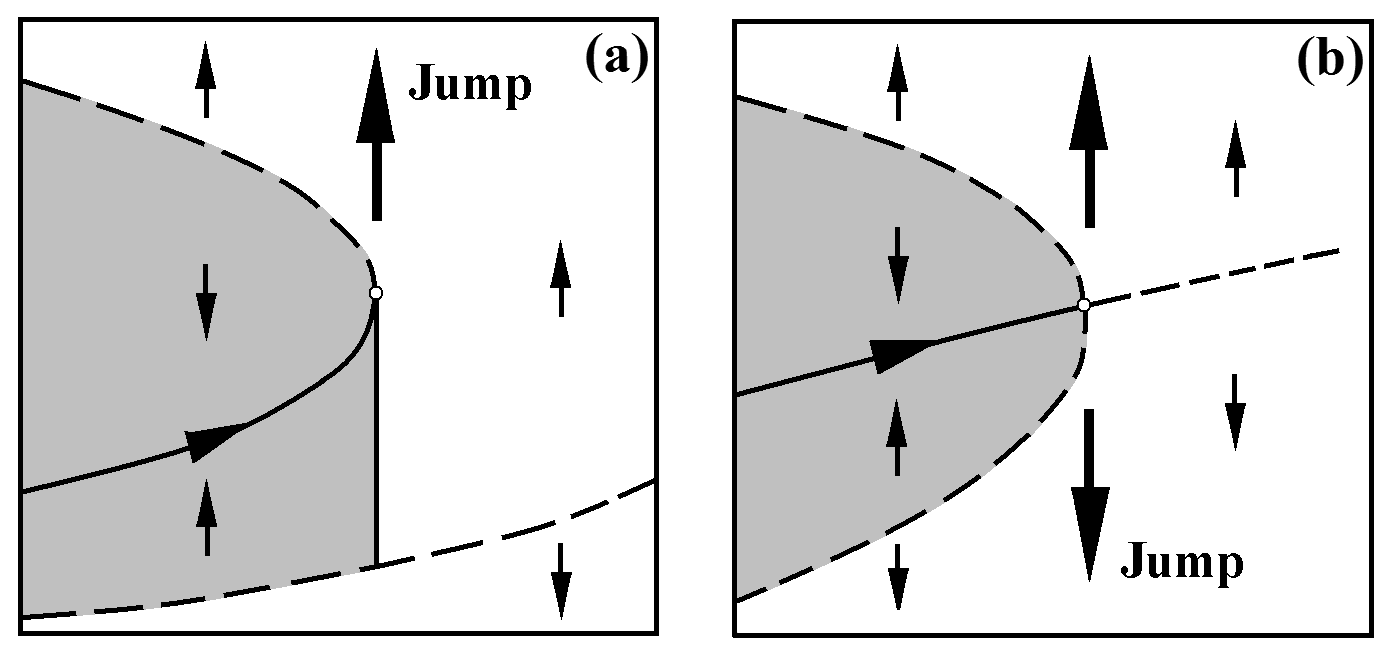}
  \caption{Basin boundary transformations at dangerous
    bifurcations. (a) the fold bifurcations and (b) the subcritical
    bifurcations. Solid curves denote stable paths while broken curves
    denote unstable paths. The basin is shown in grey.}
  \label{fig:c360}
\end{figure}

The relevant bifurcations are the so-called \emph{codimension-one}
events (\citet{Thompson1994,Thompson2002}), namely those that can be
typically encountered under a gradual change of a single control
parameter. The complete list of the dangerous codimension-one
bifurcations is given in Table~\ref{tab:bifs}, where we give a brief
description of the nature of the shrinking basin. More details can be
found in \citet{Thompson2010}.
\begin{table}
  \centering
  \begin{tabbing}
    \textbf{(a) L}\=\textbf{ocal Saddle-Node Bifurcations}\hspace*{10em}\=\\
    \> Static Fold (saddle-node of fixed point) \> one-sided basin shrinkage\\
    \>Cyclic Fold (saddle-node of cycle)\>one-sided basin shrinkage\\
    \textbf{(b) Local Subcritical Bifurcations}\\
    \>Subcritical Hopf\>complete basin shrinkage\\
    \> Subcritical Neimark-Sacker (secondary Hopf)\>complete basin shrinkage\\
    \>Subcritical Flip (period-doubling)\>complete basin shrinkage \\
    \textbf{(c) Global Bifurcations}\\
	\>Saddle Connection (homoclinic connection)\>
        outside shrinkage around cycle\\
	\> Regular-Saddle Catastrophe (boundary crisis)\>
        fractal basin shrinkage\\
	\>Chaotic-Saddle Catastrophe (boundary crisis)\>
        fractal basin shrinkage
  \end{tabbing}
  \caption{Dangerous Bifurcations and the behaviour of the basin of attraction}
  \label{tab:bifs}
\end{table}
Focusing first on the local bifurcations (headings (a) and (b) of
Table~\ref{tab:bifs}) the shrinking boundaries are illustrated
schematically in Figure~\ref{fig:c360}, with the control parameter plotted
horizontally and the response plotted vertically.
In the first picture, Figure~\ref{fig:c360}(a), we illustrate both
types of the saddle-node fold bifurcation: the static fold involving
a path of equilibrium fixed-points, and the cyclic fold involving a
trace of periodic orbits. For these the basin shrinkage is
one-sided. The basin is bounded by the unstable side of the
parabolically folding path.

Next, in Figure~\ref{fig:c360}(b), we illustrate the local subcritical
bifurcations of which the Hopf, flip and Neimark are codimension-one
(we remember that the more familiar static pitchfork bifurcations are
codimension-two, being structurally unstable against a
symmetry-breaking perturbation). Here we see an unstable solution,
which controls the basin boundary, shrinking parabolically around the
stable solution, from which the (noise-free) system jumps out of our
field of view as the control parameter, $\mu$, reaches its critical value
of $\mu_C$.

\begin{figure}[th]
  \centering
  \includegraphics[scale=0.15]{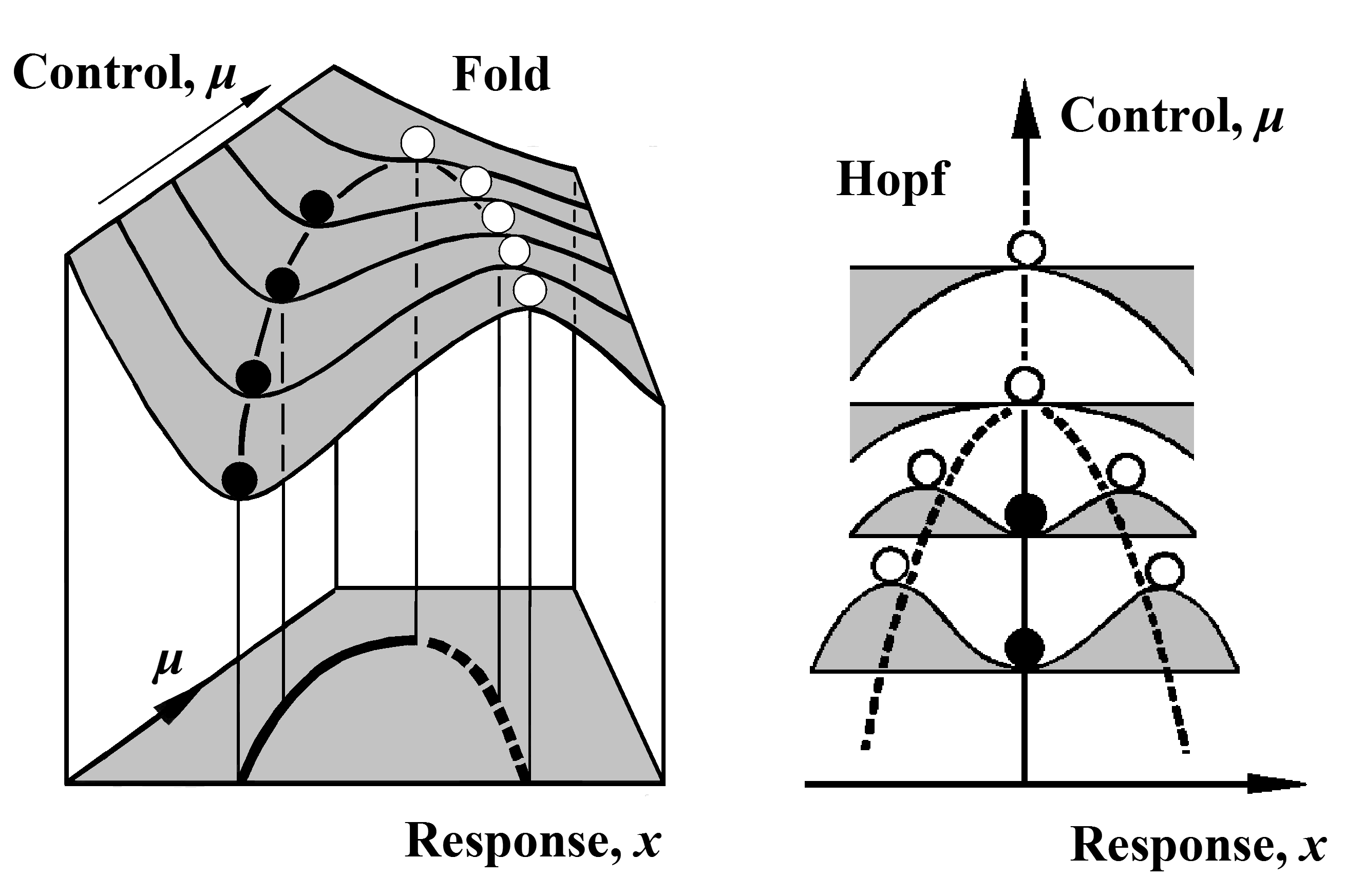}
  \caption{Total potential energy transformations in (a) the
    saddle-node fold and (b) the Hopf bifurcation (via
    averaging). Black balls denote stable equilibrium states while
    white balls denote unstable states.}
  \label{fig:c361}
\end{figure}

The static fold and the Hopf bifurcation have a theoretically
well-defined potential energy surface that neatly summarizes the
shrinking basin as illustrated in Figure~\ref{fig:c361}. Note that
Figure~\ref{fig:c361}(a) will be discussed more fully in
section~\ref{sec:concepts}(\ref{sec:fold}).
\jsrem{For other dangerous bifurcations, where the existence of an underlying
potential surface lacks any theoretical backing (and indeed may be
technically impossible because of fractal features in the dynamics)
there are nevertheless occasions when it is useful and practical to
consider such a surface for escape predictions in the presence of
significant noise as we shall demonstrate in section%~\ref{sec:spiral}.
}
\subsection{A closer look at the fold}
\label{sec:fold}
The fold is the simplest and most common way in which an equilibrium
of a nonlinear dissipative dynamical system can lose its
stability. Having only a single active degree of freedom ($x$) and
being observable under the variation of a single control parameter
($\mu$) it can be illustrated on a graph of $x$ against $\mu$ as a
smooth curve that simply reaches an extreme value (a maximum, say) of
the control $\mu$, as shown in Figure~\ref{fig:c360}(a). So although
it is traditionally called a saddle-node bifurcation, it does not
exhibit any obvious `bifurcation' of a path, but rather just a smooth
folding of an existing path. In the case of a static fold (to which we
shall largely restrict our attention) the path is a trace of
equilibrium fixed points, while for a cyclic fold the path would be a
trace of periodic solutions.

Assuming that we are close to a fold, the intrinsic damping of the
system will have become super-critical so the system will have the
non-oscillatory response of an over-damped particle sliding (in thick
oil, as we might imagine) on a notional potential energy surface
$U(x,\mu)$. This is illustrated for the fold in
Figure~\ref{fig:c361}(a), where we show the potential energy surface
erected over the $(x,\mu)$ base plane. Here the solid part of the
equilibrium path, corresponding to energy minima, is a stable node,
while the dashed part, corresponding to energy maxima (more generally
a geometrical saddle in higher dimensions) is an dynamically unstable
saddle.

We see that as the control $\mu$ is slowly increased towards its
critical value, $\mu_C$, the stable equilibrium solution gets closer
and closer to the hill-top potential barrier, and so gets
progressively more precarious in the presence of noisy disturbances.
With $\mu$ increasing at an infinitesimal rate, noise will ensure
escape over the hill-top before $\mu$ reaches $\mu_C$. If, however,
$\mu$ increases rapidly, this early escape may be delayed or supressed
as we have examined quantitatively elsewhere (\citet{Thompson2011}).
\begin{figure}[th]
  \centering
  \includegraphics[scale=0.3]{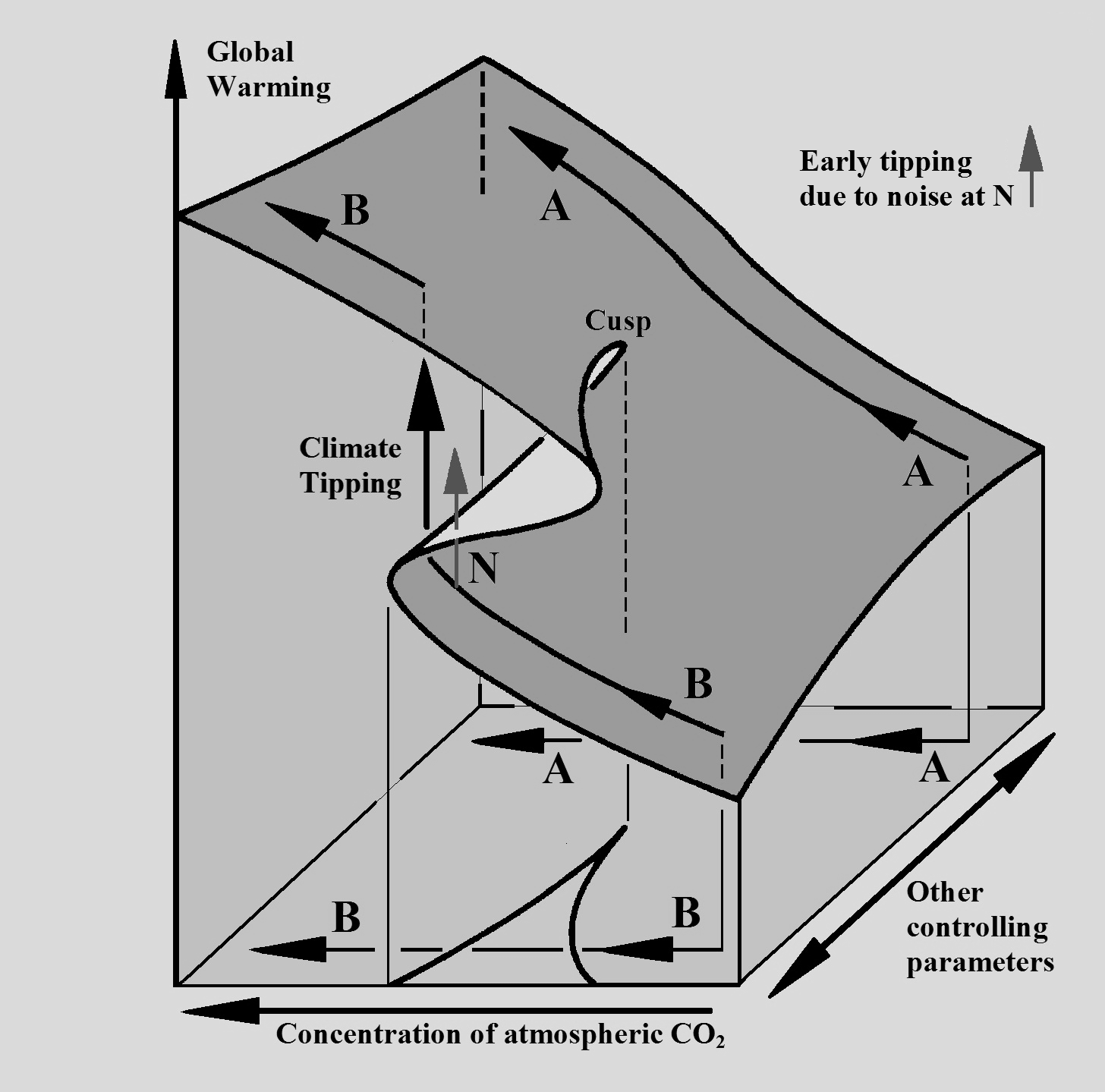}
  \caption{The cusp catastrophe and its associated pattern of
    folds. The notional axes show how such a cusp might arise in a
    climate tipping scenario.}
  \label{fig:c357}
\end{figure}

To give a greater perspective to our view of the fold it is useful to
introduce another control parameter $\mu_2$ such that our parameter
space is now represented jointly by the $(\mu_1,\mu_2)$-plane. We can
now erect an equilibrium surface, of height $x$, over the
two-dimensional control plane as illustrated in
Figure~\ref{fig:c357}. Here we see two fold lines whose projections
divide the control plane into regimes exhibiting either one or three
co-existing solutions. These fold lines coalesce and vanish at the
cusp. This cusp is a codimension-two phenomenon, which is typically
only observed when we have independent control of both parameters
($\mu_1$ and $\mu_2$). This is made transparent by the fact that a
typical path in control space cannot be expected to pass precisely
under the cusp point; by contrast, if a parameter path passes through
a fold line all near-by paths also cross this fold line, which is what
`the fold has codimension-one' means. For an appropriate choice of
path in the ($\mu_1,\mu_2)$-plane through the cusp one sees a
pitchfork bifurcation (\citet{Thompson2002}).

We will consider the case in which the main equilibrium sheet (in
grey) is stable while the narrow inverted sheet (white) is
unstable.

A climate system will have a large number of possible control
parameters, but in trying to predict a tipping point we will be
studying a particular evolution which corresponds to a particular
route in control space; this is precisely why we are only likely to
encounter a fold, and not a cusp in Figure~\ref{fig:c357}. To
explore possible scenarios that might underlie a tipping study, it is
now useful to examine some possible routes in our two-parameter model
of Figure~\ref{fig:c357}.

Consider, first, the simple scenario in which the system parameters
take path $\overline{BB}$ leading transversally through a fold
line. This is the classical encounter with a fold that is implicitly
envisaged in earlier work. In the presence of noise a premature tip
from N occurs with a certain probability depending on the noise level
and the speed with which the path is traversed, as analysed by
\citet{Thompson2010}. If the noise level were particularly high, a
jump from N could easily be perceived as a purely noise-induced jump
with no bifurcation, even though it is actually caused by the
proximity to the fold: this could hold even if there were no movement
in the control space at all, the system having rested at N throughout
(\ashwin call this N-tipping). Another scenario, not
specifically illustrated in the figure, could arise if the route in
control space approached the fold line, but then turned back away from
it. Analysis would then show a temporary decrease of the strength of
the attraction (with the AR(1) coefficient moving towards unity) which
might be discounted as a false alarm, even though the system did
approach a fold and a noise-induced escape was temporarily probable.
% In this context, see \citet{Biggs2009}.

Finally, in the scenario of path $\overline{AA}$, one would observe a
temporary increase of the AR(1) coefficient, even though no
bifurcation is apparent but rather a gradual shift of the steady
state. This shows that a rapid transition can be related to a nearby
fold that was just barely missed. Similarly, the scenario called
R-tipping by \ashwin occurs if one changes the other
control parameter in Figure~\ref{fig:c357} sufficiently rapidly from
the path $\overline{AA}$ to N (avoiding the folds). \jsmod{In this
  scenario the system jumps to the other stable equilibrium despite
  having never crossed a bifurcation curve.}

\section{Time series analysis for leading order terms --- qualitative
  overview}
\label{sec:timeseries}
The extraction of the underlying dynamics from time series is divided
into subtasks of increasing difficulty. Assume that we have a time
series coming from a process that is either stationary or has a slowly
drifting underlying parameter $\mu$ as suggested by the paths in
Figure~\ref{fig:c357}. In order to predict (or identify) the dynamics
one attempts to extract the leading orders of the dominant components
of the underlying equations of motion.
\jsmod{
\begin{figure}[tbh]
  \centering
  \includegraphics[width=0.9\textwidth]{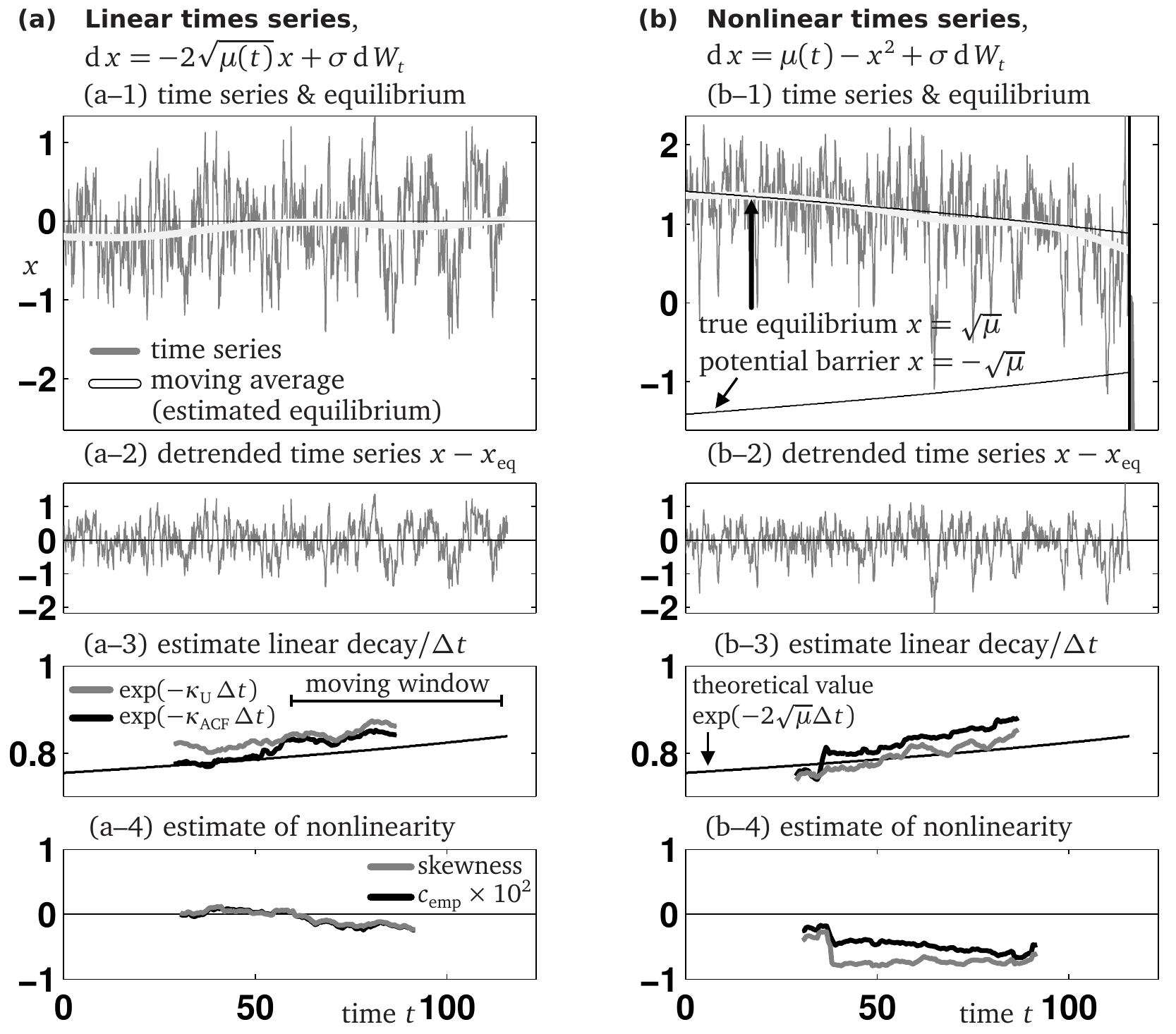}
  \caption{Comparison of time series generated by linear (left
      column, (a-1) to (a-4)) and nonlinear process (right column,
      (b-1) to (b-4)). Both time series are stationary and have
      identical linear decay rates. As measures for nonlinearity we
      used skewness and the coefficient $c_\mathrm{emp}$ (as defined
      below in equation~\eqref{eq:fp1}). Parameters: $\mu(0)=2$,
      $N=1\,223$, $\Delta t=0.1$, $\sigma=1$, window length for linear
      and nonlinear analysis, $w=N/2$, $\epsilon=0.01$.}
  \label{fig:snfdata}
\end{figure}
} We illustrate the steps using \jsmod{a time series} obtained from a
process which has \jsmod{the saddle-node normal form with a slowly
  drifting normal form parameter $\mu$} as its deterministic part
(corresponding to Figure~\ref{fig:c361}(a)) and is influenced by
additive Gaussian noise. We permit the parameter $\mu$ to drift slowly
with speed $\epsilon$:
\begin{equation}
  \label{eq:snfnoise}
  \begin{split}
    \d x&=\left[\mu-x^2\right]\d t +\sigma \d W_t \mbox{, where}\\
    \d\mu& =-\epsilon \d t\mbox{.}
  \end{split}
\end{equation}
Here $x$ is the response, $\mu$ is the control parameter, and $t$ is
the time. The Gaussian perturbation $\d W_t$ has zero mean and unit
variance such that $\sigma$ controls the noise induced variance added
to the dynamics of the deterministic saddle-node form $\dot
x=\mu-x^2$, which has the stable equilibrium
$x_\mathrm{eq}=\sqrt{\smash[b]{\mu}}$. The second equation prescribes
the linear sweep of $\mu$ from its starting value $\mu_0>0$ towards
zero at rate $\epsilon\ll1$. For fixed $\mu>0$ the unperturbed system
($\sigma=0$) has a stable equilibrium at
$x=x_\mathrm{eq}=\sqrt{\smash[b]{\mu}}$ and an unstable equilibrium at
$x=-\sqrt{\smash[b]{\mu}}$.

In Section~\ref{sec:quant} we will use the saddle-node normal form
\eqref{eq:snfnoise} to \jsmod{test the ability of a range of estimators to
distinguish a time series generated by a process with a quadratic
nonlinearity from a linear time series. Two example time series
  and their analysis are shown in Figure~\ref{fig:snfdata}. Both time
  series have identical linear properties but differ in the dominant
  nonlinear term of the deterministic part of the dynamics.The
  softening of the full (estimated) nonlinear potential well
  (corresponding to the illustration in Figure~\ref{fig:c361}(a)) is
  displayed in Figure~\ref{fig:c365}.  }

\jsrem{The left column (a) of Figure~\ref{fig:snfdata} shows a linear
  time series generated by \eqref{eq:snfnoise} with densely taken
  measurements (length of time series $N=28,803$, time step $\Delta
  t=0.1$, drift speed $\epsilon=0.01$ for system parameter $\mu$). The
  right column shows the results for a short time series (also
  generated by \eqref{eq:snfnoise}) with sparser measurements.}
\subsection{Equilibrium --- zero order}
\label{sec:0order}
\noindent The location of the slowly drifting equilibrium and its
trend is estimated first. \lenton call this step detrending
and propose either piecewise linear fitting or filtering with a
Gaussian kernel. The result of filtering with a Gaussian kernel is
shown as a thick grey curve in the row 1 of
Figure~\ref{fig:snfdata}. Every method for detrending requires a
bandwidth parameter (for example the width of the Gaussian
kernel). For Figure~\ref{fig:snfdata} we chose the bandwidth which the
kernel density estimation routine (Matlab's \texttt{kde} by
\citet{Botev2010}) offered. The thin black curve in
\jsmod{Figure~\ref{fig:snfdata}(b-1)} shows the true value of the
quasi-equilibrium (for $\epsilon=0$). It gives a visual estimate of how the
reliability of the zero-order estimates depends on the quality of the
data. The oscillatory deviations of the extracted equilibrium (thick
light grey curve) from its true value (thin black curve) shows that
the automatically chosen bandwidth was slightly smaller than the
theoretical optimum.

\subsection{Linear decay rate --- first order}
\label{sec:1order}
\noindent Assuming that the deterministic part of the underlying
process has a stable equilibrium $x_\mathrm{eq}$ (which is possibly
slowly drifting), and that the zero-order estimate is an approximation
for $x_\mathrm{eq}$, the next step is to estimate the linear decay
rate $\kappa$ toward $x_\mathrm{eq}$ and the trend of $\kappa$ over
time. Generally, if the underlying deterministic process is
high-dimensional then these estimates capture the dominant (that is,
smallest) decay rate. They are typically applied to the detrended time
series, that is, $\tilde x=x-x_\mathrm{eq}$, which is shown in the
second row of Figure~\ref{fig:snfdata}.
\begin{itemize}
\item \textbf{ACF} The $k$-step (usually one-step) autocorrelation
  function (ACF) $\alpha$ is fitted directly to the detrended time
  series: $\tilde x_{k+1}=\alpha \tilde x_k$. This was introduced by
  \citet{Held2004} as degenerate finger-printing for the analysis of
  climate time series. The ACF $\alpha$ is related to $\kappa$ via
  \begin{displaymath}
    \alpha=\exp(-\kappa_\mathrm{ACF}\,\Delta t)\mbox{.}
  \end{displaymath}
  The estimate for $\kappa_\mathrm{ACF}$ is shown in
  Figure~\ref{fig:snfdata}(a-3) and (b-3). 
\item \textbf{DFA} Detrended fluctuation analysis has been introduced
  by \citet{Livina2007} for climate time series, see also
  \lenton for a short description of how one extracts $\kappa$
  from this analysis.
\item \textbf{Quasi-stationary density} If one assumes that the
  deterministic dynamics of the detrended time series $\tilde x$ is
  essentially that of an overdamped particle moving in a potential
  well $U(\tilde x)$, that is,
  \begin{equation}
    \label{eq:potwell}
    \d \tilde x=-\partial_xU(\tilde x)\d t+\sigma \d W_t\mbox{,}
  \end{equation}
  then the linear decay rate equals the second derivative
  $\partial_{xx}U(0)$ of $U$ at the bottom of the well, $\tilde
  x=0$. Moreover, the stationary density $p(x)$ is related to the
  shape of the potential well $U$ via the stationary Fokker-Planck equation
  \begin{align}
    \frac{1}{2}\partial_{xx} p(\tilde x)&=
    \partial_x\left[-\sigma^{-2}\partial_xU(\tilde x)\,p(\tilde
      x)\right]\mbox{, which gives, integrated once,}\nonumber\\
    \label{eq:fp}
    \frac{1}{2}\partial_x p(\tilde x)&=
    -\sigma^{-2}\partial_xU(\tilde x)\,p(\tilde x)+\sigma^{-2}c\mbox{,}
  \end{align}
  where $c$ is the constant of integration. This constant satisfies
  $c=\lim_{\tilde x\to\pm\infty}\partial_xU(\tilde x)p(\tilde x)$, and
  can be interpreted as the flow rate ($c<0$ indicates flow toward
  $-\infty$). Remembering that the drift speed $\epsilon$ of the
  system parameter $\mu$ is small, the true time-dependent probability
  density $p$ still satisfies \eqref{eq:fp} approximately at each
  $\mu$:
  \begin{equation}
    \label{eq:fpdrift}
    \frac{1}{2}\partial_x p(\tilde x,\mu)=
    -\sigma^{-2}\partial_xU(\tilde x,\mu)\,p(\tilde x,\mu)+
    \sigma^{-2}c(\mu)+O(\epsilon)\mbox{.}    
  \end{equation}
  In \eqref{eq:fpdrift} we included the dependence of all quantities
  on $\mu$ expressly and remember that $\dot \mu(t)$ is of order
  $\epsilon$. Assuming quasi-stationarity, we neglect the
  order-$\epsilon$ terms and determine the coefficients,
  $-\sigma^{-2}\partial_xU$ and $\sigma^{-2}c$, by fitting the
  empirical density $p_\mathrm{emp}(x)$ from a time window
  $[t-w/2,t+w/2]$ to the Equation~\eqref{eq:fpdrift}. Specifically,
  $c_\mathrm{emp}=\sigma^{-2}c$ is a scalar and $\partial_xU(0)$ is
  equal to $0$ after detrending such that one has to fit two
  coefficients, $\kappa_U$ and $c_\mathrm{emp}$ using
  \eqref{eq:fpdrift} if one truncates $\partial_xU$ after first order:
  \begin{equation}
    \label{eq:fp1}
    \frac{1}{2}\partial_x p_\mathrm{emp}(\tilde x)=
    -\kappa_U \tilde x\,p_\mathrm{emp}(\tilde x)+c_\mathrm{emp}\mbox{.}
  \end{equation}
  In problems where no escape is possible one can drop the term $c$ in
  \eqref{eq:fp} (and, thus, in \eqref{eq:fpdrift} and \eqref{eq:fp1}),
  and simplify \eqref{eq:fp} to
  \begin{equation}\label{eq:logu}
    U=-\frac{\sigma^2}{2}\log p\mbox{.}
  \end{equation}
  \citet{Livina2011} used relation \eqref{eq:logu} (fitting of $U$
  with higher-order polynomials) to detect potential wells in time
  series that included rapid transitions. Fitting $U$ with a parabola
  (where $U(0)$ is irrelevant and $\partial_xU(0)=0$) is equivalent to
  using the variance $\sigma^2_\mathrm{emp}$ of $p_\mathrm{emp}(\tilde
  x)$ if the density $p_\mathrm{emp}$ is
  Gaussian. \citet{Ditlevsen2010} state that monitoring the empirical
  variance $\sigma^2_\mathrm{emp}$ helps to avoid erroneous detection
  of false trends.
\end{itemize}
\lenton compare these three estimates,
\jsmod{$\kappa_\mathrm{ACF}$, the DFA-based estimate for $\kappa$, and
  the variance $\sigma^2_\mathrm{emp}$ (which is inversely
  proportional to $\kappa_U$),} in detail using time series arising in
climate models and in palaeoclimate records. Methods based on
properties of the spectrum of the time series were proposed and
investigated by \citet{Kleinen2003,Biggs2009} but are not discussed
here.  \jsmod{Figure~\ref{fig:snfdata} shows in row~3 how the
  estimates of $\kappa_U$ from the quasi-stationary density compare to
  the ACF estimates $\kappa_\mathrm{ACF}$. For both time series the
  estimates are quantitatively close to the theoretical value of
  $\kappa$ (which is indicated by the thin black line in
  Figure~\ref{fig:snfdata}(a-3,b-3)).}
\jsrem{Figure~\ref{fig:snfdata}(b-3) shows that even for sparse and
  short time series the estimates for the linear decay rate are
  qualitatively correct but that discerning trends such as in
  Figure~\ref{fig:snfdata}(a-3) is possible only for longer time
  series.}  \jsmod{We also observe that in
  Figure~\ref{fig:snfdata}(a-3,b-3) the local trends of
  $\kappa_\mathrm{ACF}$ from the autocorrelation estimate and of
  $\kappa_U$ based on the quasi-stationary density are strongly
  correlated, so it is unclear if monitoring both quantities really
  prevents false alarms.}  \jsrem{As the underlying decay rate has no
  downward trend one might call the local downward trend of the
  estimates a `false alarm'.}  \jsrem{However, the parameters of the
  system make a noise-induced escape (or N-tipping in the notation of
  \mbox{\ashwin}) likely. Escape occurs indeed at time
  $t=40$ in Figure~\ref{fig:snfdata}(b-1).}

There is one notable difference between the AR(1) and the DFA estimate
on one hand and the distribution-based estimate $\kappa_U$ on the
other hand: the estimate for $\kappa_U$ based on the quasi-stationary
density only looks at the distribution of the values $\tilde x_k$ in
the window of interest but does not care about the order in which they
appear. This is in contrast to, for example, the AR(1) estimate which
determines the correlation between each value $\tilde x_k$ and its
predecessor.  

\subsubsection*{Estimate of the input noise amplitude}
\noindent If the AR analysis shows the presence of a single positive
dominant coefficient (which gives evidence that the underlying
dynamics has really one distinct direction in which the decay is
slowest) then one can combine both estimates of the decay rate to
obtain an estimate of the amplitude $\sigma$ of the additive
noise. \jsmod{Equations \eqref{eq:fp} and \eqref{eq:fp1} imply that
  $\kappa_U$ is an estimate for $\sigma^{-2}\kappa$ where $\kappa$ and
  $\sigma^2$ are the true linear decay rate and input noise
  amplitude. Thus, if we replace the unknown $\kappa$ by the estimate
  $\kappa_\mathrm{ACF}$ we can use the relation between $\kappa_U$ and
  the true $\kappa$ to estimate $\sigma^2$:}
\begin{equation}
  \sigma_\mathrm{emp}^2=\kappa_\mathrm{ACF}/\kappa_U\label{eq:sigmaemp}
\end{equation}
where $\kappa_\mathrm{ACF}$ is the estimate of the linear decay rate
obtained from the autocorrelation and $\kappa_U$is the estimate
obtained from the quasi-stationary density $p_\mathrm{emp}$. \jsmod{This is
an alternative to using the residual of the least-squares fit that one
obtains when estimating $\kappa_\mathrm{ACF}$. We found that the
residual systematically underestimates the noise level in the normal
form examples shown in Fig.~\ref{fig:snfdata}.}

\subsection{Nonlinearity --- basin boundary}
\label{sec:nonlin}
\noindent One problem left open by the linear analysis is that the
quantitative value of the estimated linear decay rate $\kappa$ and the
estimated noise-level do not give any indication for the probability
of tipping. What the linear methods estimate is the $\alpha$ and the
$\sigma$ of a discrete process
\begin{equation}\label{eq:lindisc}
  \tilde x_{k+1}=\alpha \tilde x_k+\sigma\eta_k
\end{equation}
(where the $\eta_k$ are assumed to be independent and normally
distributed random numbers). The discrete process \eqref{eq:lindisc}
is the linearisation of the time-$\Delta t$ map of the continuous
process \eqref{eq:potwell}. An estimate for $\alpha$ less than but
close to unity does not necessarily indicate that we are close to a
stability boundary for the nonlinear problem \eqref{eq:potwell}.
Rather, it implies that the time step $\Delta t$ between successive
measurements has been small compared to the mean decay time to half,
$\log 2/\kappa$, of the process:
\begin{displaymath}
  \alpha(t_k)=\exp\left(-\kappa(t_k)\Delta t\right)=
  1-\kappa(t_k)\Delta t+O\left((\Delta t)^2\right)\mbox{.}
\end{displaymath}
The \emph{trend} of the estimated $\alpha$ is an indicator for
incipient tipping because (if extrapolated in any way) it gives an
estimate for the time to tipping if one ignores the possibility of
early escape. However, the trend of $\alpha$ is less certain even in
the artificial time series of Figure~\ref{fig:snfdata}. Similarly, the
estimated noise amplitude $\sigma$ has to be compared to the
coefficient in front of the leading nonlinear term of the right-hand
side (which equals $-1$ in \eqref{eq:snfnoise}). The ratio between
$\sigma$ and nonlinear term enters the estimates for the probability
of tipping if there is no discernible upward trend in $\alpha$
(N-tipping), or for the probability of early escape if the trend of
$\alpha$ points upwards (see \citet{Thompson2011}). So without at
least an order-of-magnitude estimate of the leading nonlinear term the
linear methods only provide an estimate for the tipping time if there
is a clearly discernible trend in $\alpha$, and even then this
estimate discounts early escape.

\jsmod{Figure~\ref{fig:snfdata} illustrates this problem. Apart from a
  slow trend both time series in Figure~\ref{fig:snfdata} have
  identical properties at the linear level. However, while for the
  time series in Figure~\ref{fig:snfdata}(b-1) the median time to
  tipping is $t=100$ according to \citet{Thompson2011} (indeed, it
  escapes shortly after $t=120$), for the time series in
  Figure~\ref{fig:snfdata}(a-1) the extrapolation (correctly) predicts
  tipping (or rather, linear instability) at $t\approx320$.  Thus, we
  observe that the linear properties of the time series alone are not
  sufficient to provide estimates about tipping, for example, the
  median time to tipping. The difference between the time series in
  Figure~\ref{fig:snfdata}(a-1) and Figure~\ref{fig:snfdata}(b-1) is
  in the dominant nonlinear term of the deterministic part. That this
  term plays a role is not surprising because tipping is a feature of
  nonlinear systems. Extracting its presence from a time series is a
  necessary step that has to follow the linear analysis.}  \jsmod{
\begin{figure}[th]
  \centering
  \includegraphics[width=0.6\textwidth]{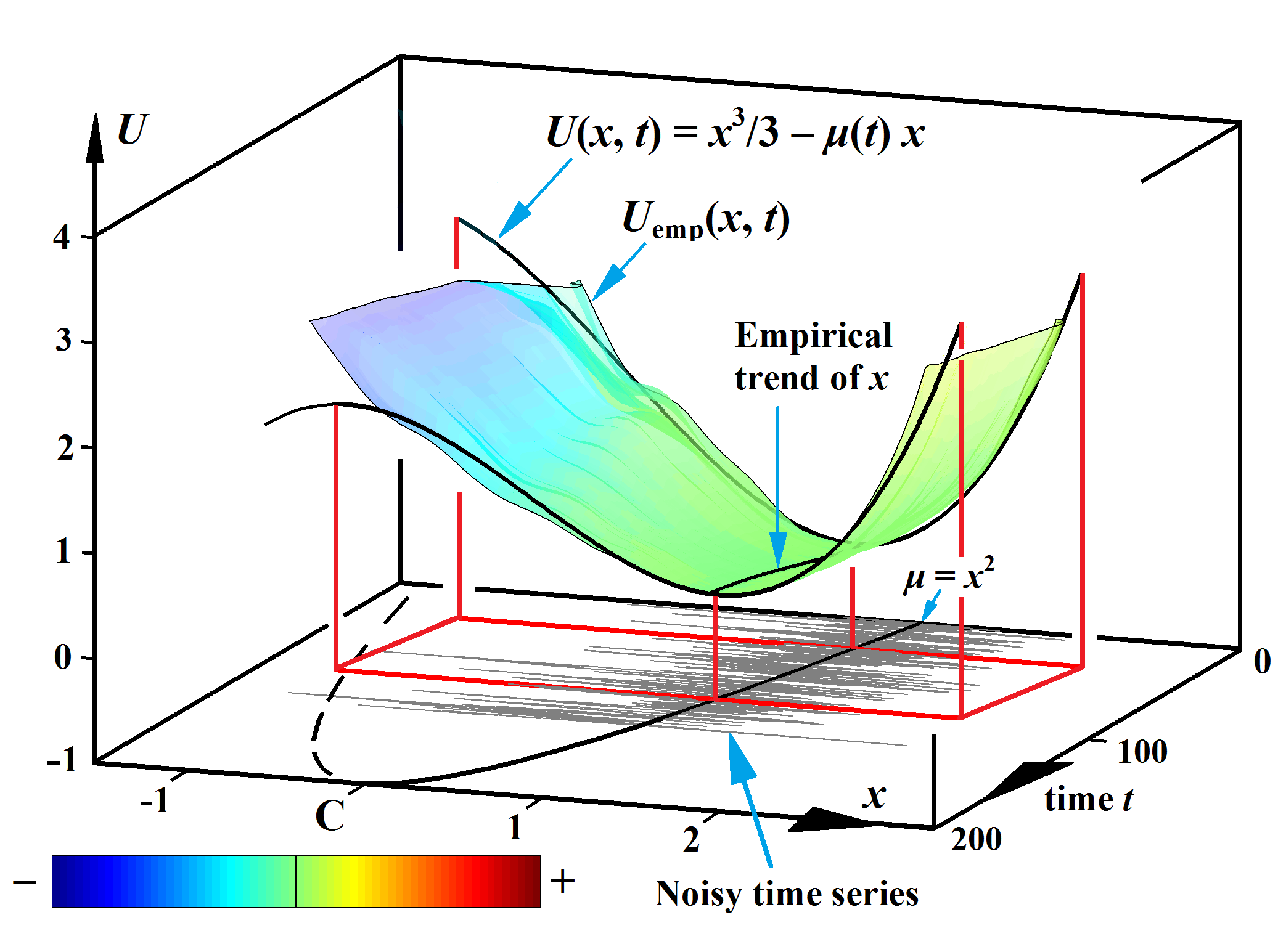}
  \caption{Nonlinear potential energy surfaces extracted from noisy
    time series using sliding windows ($w=N/2$) for the saddle node
    fold \protect\jsmod{($U_\mathrm{emp}=
      -\sigma_\mathrm{emp}^2\log(p_\mathrm{emp})/2$ where
      $\sigma_\mathrm{emp}$ was estimated using \eqref{eq:sigmaemp})}.
    At either end of the surface we have plotted the known function
    $U(x,t)=x^3/3-\mu(t)x$ ($\sigma=1$, $\epsilon=0.01$ and
    $\mu(0)=2$).  Red denotes a positive deviation, blue a negative
    deviation from the parabola given by linear theory.}
  \label{fig:c365}
\end{figure}}
The approaches presented in the following all generalise the estimate
for $\kappa_U$, based on the quasi-stationary density, by looking at
the full nonlinear potential well $U_\mathrm{emp}(\tilde x)$
associated to the quasi-stationary density $p_\mathrm{emp}(\tilde x)$
via the Fokker-Planck
equation~\eqref{eq:fpdrift}. Figure~\ref{fig:c365} illustrates how
$U_\mathrm{emp}=-\sigma^2\log(p_\mathrm{emp})/2$ looks like for the
time series shown in Figure~\ref{fig:snfdata} column (b). The
realisation of the noisy, evolving time series is shown in the base
plane together with the equilibrium path $\mu=x^2$, and the estimated
(empirical) potential, $U_\mathrm{emp}$ is shown as the coloured
surface. The estimate for the potential well is taken in a sliding
window, which explains why the surface can only be determined in the
central region of the time series with half the length of the window
unrepresented at either end. The numerically estimated equilibrium
trend $x_\mathrm{eq}$ is displayed as a black curve at the bottom of
the valley of $U_\mathrm{emp}$. For comparison with this
$U_\mathrm{emp}$ we display at either end the real (nonlinear)
potential $U(x,t)$ corresponding to \eqref{eq:snfnoise}, which at the
later time is already showing the hill-top to the left. % The fit
% between $U_\mathrm{emp}$ and $U(x,t)$ is seen to be very
% satisfactory indeed.
The colour in Figure~\ref{fig:c365} shows the deviation of the
empirical potential from the linear-theory parabola, which is not
itself displayed. The softening, highlighted by the blue colouration,
is seen on the left hand side which is approaching the
hill-top. Transparency of the colour is used to show the number of
data points that support a particular part of the surface
$U_\mathrm{emp}(\tilde x)$ (as given by the empirical quasi-stationary density
$p_\mathrm{emp}(\tilde x)$).

Several methods have been proposed to detect nonlinearity in the
potential well in the form of this type of softening (or tilting).
\begin{itemize}
\item \textbf{Skewness} \citet{Guttal2008,Guttal2008a} studied how the
  closeness of a bifurcation influences the skewness $\gamma$ of the
  empirical quasi-stationary density from time series. As this
  approach is also based on the empirical quasi-stationary density
  $p_\mathrm{emp}$, it generalises the estimate $\kappa_U$ to
  non-parabolic features of the well $U$ and non-Gaussian features of
  the stationary density. \citet{Guttal2008,Guttal2008a} proposed to
  look at the \emph{trend} of the skewness $\gamma$ over time to
  detect incipient bifurcations. \jsmod{Figure~\ref{fig:snfdata}(b-4)
    shows how the skewness, taken over a single window, changes. It is
    noticeably negative (compare to the empirical skewness from the
    linear time series in Fig~\ref{fig:snfdata}(a-4)) but no trend is
    discernable.}
\item \textbf{Quasi-stationary density} \citet{Livina2011} generalised
  the potential well analysis for the quasi-stationary density beyond
  the estimate $\kappa_U$ for the linear decay rate, fitting the
  potential well to higher-order polynomials of even degree using
  relation~\eqref{eq:logu}. This analysis was not used to attempt
  prediction of transition probabilities from time series but it was
  applied to time series that already included all transitions to
  count the number of wells of $U$ depending on time $t$ (which
  corresponds to the number of modes (peaks) of the empirical
  quasi-stationary density $p_\mathrm{emp}$). See also
  \beaulieu for methods that detect change points in time
  series.
\item \textbf{Drift ratio} If one assumes that the underlying system
  parameter is approaching a fold more or less linearly
  ($\dot\mu(t)=\epsilon>0$) and one is sufficiently close to the fold
  already then the ratio between the drift $\Delta x_\mathrm{eq}$ of
  the equilibrium estimate and the increase of the estimated linear
  decay rate $\Delta \kappa$ estimates the quadratic term in the
  normal form \citep{Thompson2011}.
\end{itemize}
% The general area of nonlinear time series analysis developed in detail
% by \citet{Kantz2003} is more concerned with the problem of
% distinguishing a chaotic trajectory from a trajectory close to an
% equilibrium disturbed by noise, so most of its methods are not
% directly applicable to detection of nearby basin boundaries in
% noise-dominated time series.

\section{Quantitative estimates of the nonlinear coefficients}
\label{sec:quant}
The colour shading of Figure~\ref{fig:c365} suggests that the
nonlinear term in the underlying deterministic system shows up in the
quasi-stationary density $p_\mathrm{emp}$. This raises two questions.
First, can the level of deviation from linear theory be distinguished
with confidence from random chance? That is, can one state, when
looking at a time series, that a significant quadratic term is
present?  Second, is it possible to quantify the size of the nonlinear
term with reasonable level of certainty?

Finding a significant nonlinearity is far easier than estimating its
precise value because biases introduced during the analysis of the
equilibrium and decay rate (zero-order and first-order analysis) may
distort the value of the nonlinear estimate but may still keep it
significantly different from zero. We will demonstrate this in
Section~\ref{sec:geolog} for two paleo-climate time series.

\jsmod{\subsection{Properties of the stationary density of the
    saddle-node normal form}}
\label{sec:fokkerplanck}
\begin{figure}[th]
  \centering
  \includegraphics[width=0.8\textwidth]{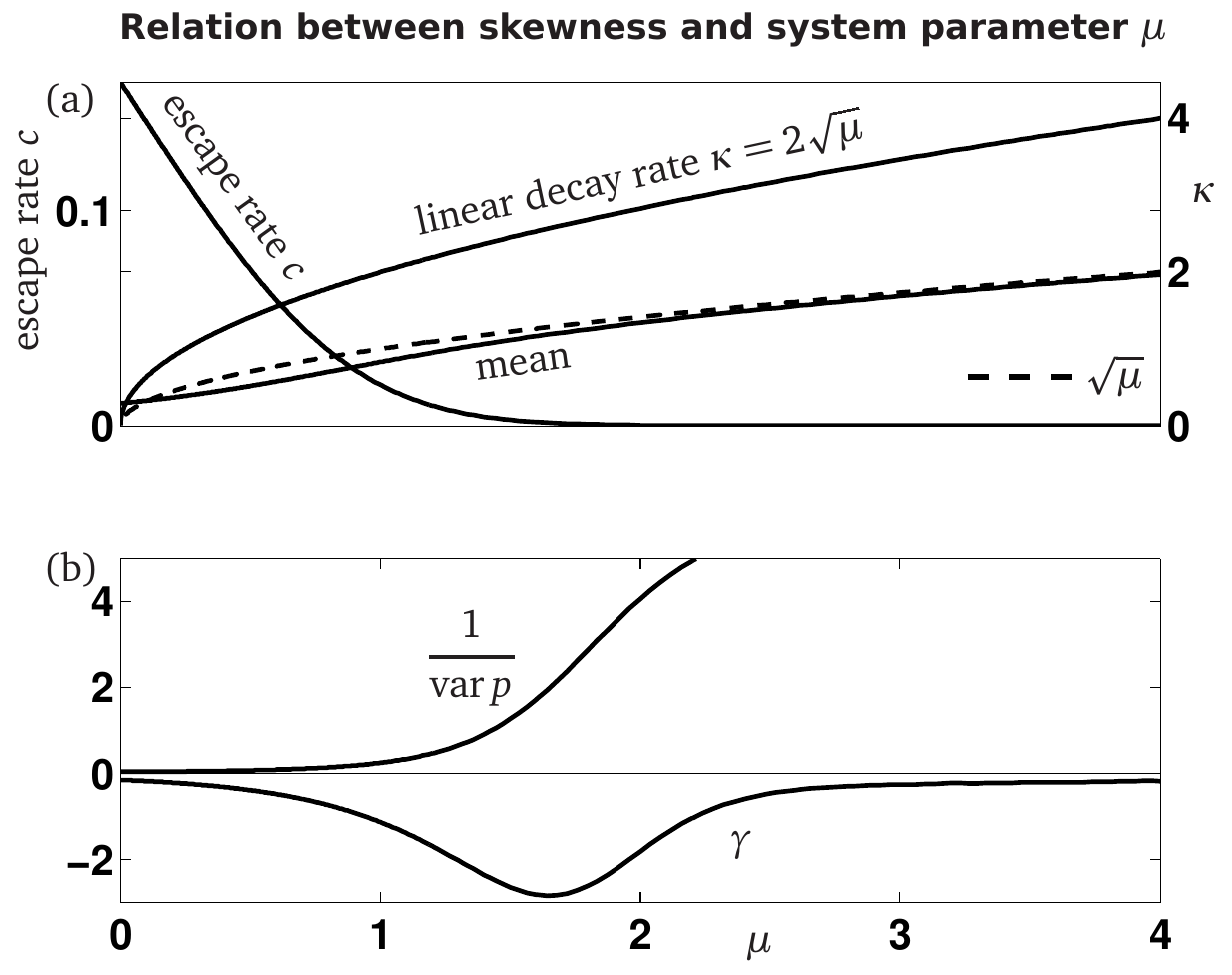}
  \caption{Dependence of mean, variance $\var p$ and skewness $\gamma$
    of the empirical stationary density for the saddle-node normal
    form with additive Gaussian noise \eqref{eq:snfnoise} \protect\jsmod{as
      computed directly from the Fokker-Planck equation~\eqref{eq:fp}
      with $-\partial_xU(x)=\mu-x^2$ and $\sigma=1$. The escape rate
      $c$ is defined by the constant of integration in
      Equation~\eqref{eq:fp}.}}
  \label{fig:skewness}
\end{figure}
Figure~\ref{fig:skewness} shows the dependence of the mean, the
variance $\var p$, and the skewness $\gamma$ of the empirical
stationary density $p_\mathrm{emp}$ on the system parameter $\mu$ for
\eqref{eq:snfnoise} as computed using the stationary Fokker-Planck
equation \eqref{eq:fp}. Note that this % is the stationary density for
% the process where the density is re-scaled to have a unit integral
% whenever a particle escapes to $-\infty$.
\jsmod{corresponds to the stationary process where the
  random realisation is re-injected at $+\infty$ whenever it has
  escaped to $-\infty$.} % Thus, the results shown in
% Figure~\ref{fig:skewness} refer to the density of realisations that
% have not escaped to $-\infty$.
We observe that the dependence of
$\gamma$ is strongly nonlinear for the saddle-node normal form with
noise ($\dot\mu(t)=\epsilon=0$, $\sigma=1$ in
Equation~\eqref{eq:snfnoise}).  This means that trends of the skewness
are likely to be difficult to ascertain when approaching a
fold. However, the presence of skewness is an indicator for the
nonlinearity uniformly for all parameters $\mu$, and its sign
indicates the direction of escape. 

Note also how the nonlinearity affects the empirical mean and the
variance $\var p$, which is no longer inversely proportional to the
linear decay rate $\kappa$. This shift of the mean and the additional
variance is in part due to the fraction of trajectories that escape,
in part due to the non-parabolic shape of the well at some distance
from its local minimum.  

\jsmod{\subsection{Practical estimates of the nonlinear term from time
    series}}
\label{sec:estnonlin}
We can estimate the nonlinearity directly as the deviation of the
empirical potential well $U$ from a parabola as suggested by
figure~\ref{fig:c365} (we are only interested in $\partial_xU$). The
simplest approach (apart from looking at skewness $\gamma$) is to fit
$\kappa_U$ and $c_\mathrm{emp}$ to the empirical density
$p_\mathrm{emp}$ using Equation~\eqref{eq:fp1}. If $\partial_xU(\tilde
x)$ is indeed linear then $c_\mathrm{emp}$ would be zero such that
$c_\mathrm{emp}$ is a signed scalar measure for the deviation of
$\partial_xU(\tilde x)$ from linearity, serving as a proxy for the
present nonlinearity. A second alternative is to expand
$-\partial_xU(\tilde x)$ to second order, incorporating a quadratic
nonlinearity directly into $\partial_xU$ with an unknown coefficient
$N_2$. This leads to
\begin{align}
  \nonumber -\sigma^{-2}\partial_xU(\tilde x)&= -\kappa_U \tilde
  x+N_2\tilde x^2\mbox{,\quad such that we fit}\\
  \label{eq:fp2}
  \frac{1}{2}\partial_x p_\mathrm{emp}(\tilde x)&=
  \left[-\kappa_U \tilde x+N_2\tilde x^2\right]\,
  p_\mathrm{emp}(\tilde x)+c_{\mathrm{emp},2}\mbox{.}
\end{align}
Both quantities, $c_\mathrm{emp}$ and $N_2$, measure the deviation of
$-\sigma^{-2}\partial_xU(\tilde x)$ from its linear approximation
$\kappa_U\tilde x$. The main difference between them is the weighting
($\tilde x^2$ for $N_2$, uniform for $c_\mathrm{emp}$). \jsmod{In
  summary, the step-by-step procedure to extract an estimate of the
  quadratic term (or a proxy for it) from a time series $x_k$ is the
  following.
  \begin{enumerate}
  \item Detrend the time series $x_k$ (see
    Section~\ref{sec:timeseries}(\ref{sec:0order})). We call the
    detrended time series $\tilde x_k$. The following steps are
    performed for each suitable time $t$ to obtain $\kappa_U(t)$,
    $c_\mathrm{emp}(t)$, $N_2(t)$ or skewness $\gamma(t)$:
  \item Choose a length $w$ of sliding windows and obtain the
    empirical stationary density $p_\mathrm{emp}(\tilde x)$ for the
    sliding window centered at $t$. (We used \texttt{kde1d} to obtain
    $p_\mathrm{emp}(\tilde x)$. For non-stationary time series we used
    the sliding window length $w=N/2$.)
  \item Use relation \eqref{eq:fp1} to find $\kappa_U(t)$ and
    $c_\mathrm{emp}(t)$ from $p_\mathrm{emp}$. Use relation
    \eqref{eq:fp2} to find $N_2(t)$ (and another estimate for
    $\kappa_U$) from $p_\mathrm{emp}$. Compute the empirical skewness
    $\gamma(t)$ from $p_\mathrm{emp}$.
  \end{enumerate}
}

\jsmod{
\subsection{Comparison of estimates using saddle-node normal form}
\label{sec:snftest}
}
\begin{figure}[th]
  \centering
  \includegraphics[width=0.9\textwidth]{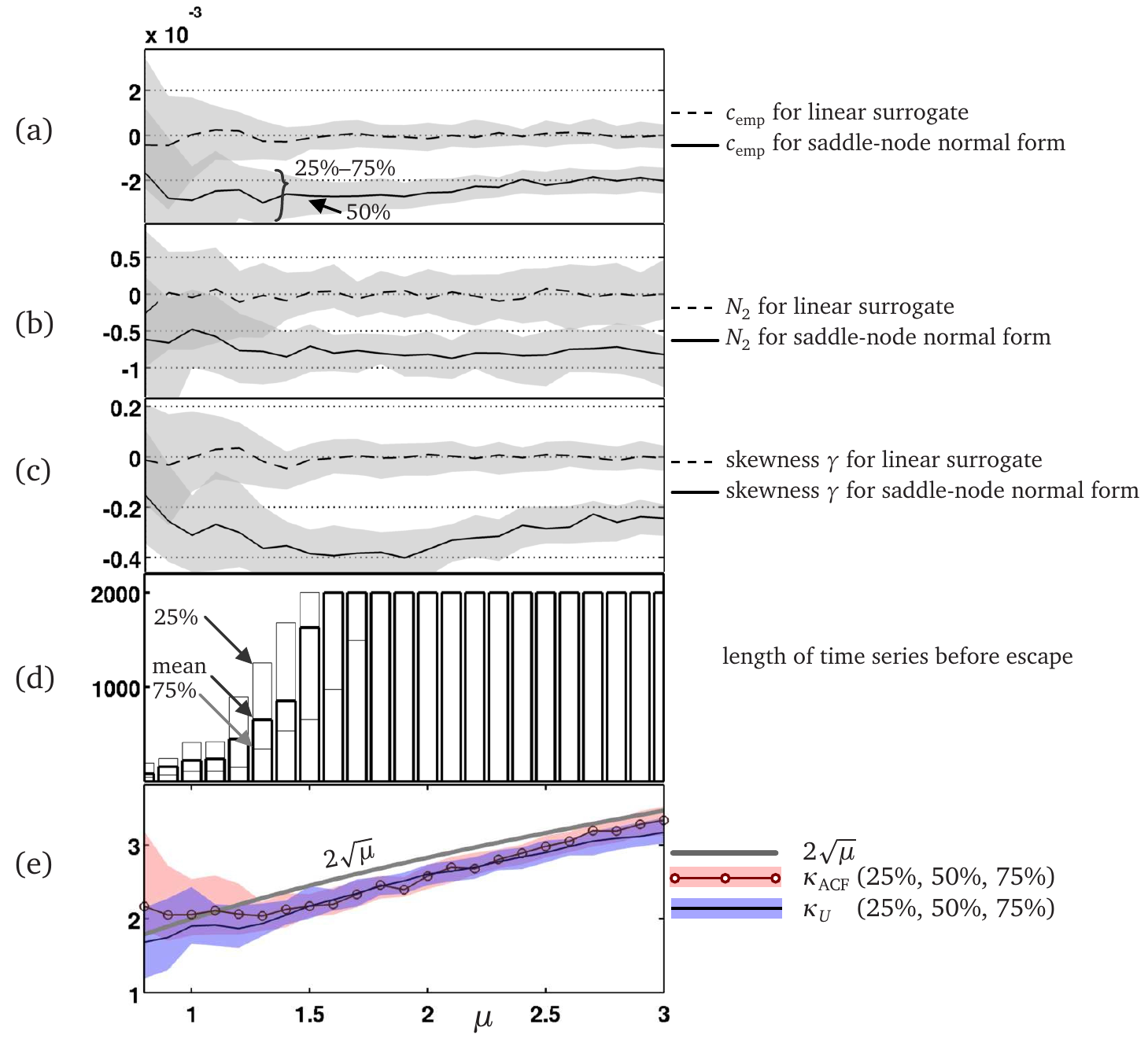}  
  \caption{Quartiles of lengths (d), and estimates $c_\mathrm{emp}$
    (a), $N_2$ (b), skewness $\gamma$ (c), $\kappa_\mathrm{ACF}$, and
    $\kappa_U$ (both (e)) for time series generated by the saddle-node
    normal form \eqref{eq:snfnoise} with noise ($\epsilon=0$, $\mu=1$)
    and for linear time series generated by
    \eqref{eq:ou}. \protect\jsmod{Without nonlinearity transition would occur
      at $\mu=0$.}}
  \label{fig:snftest}
\end{figure}
Figure~\ref{fig:snftest} shows how the estimates
for nonlinearity behave for the saddle-node normal form with noise,
Equation~\eqref{eq:snfnoise} with $\epsilon=0$ and
$\sigma=1$. % It shows $c_\mathrm{emp}$ in panel (a), $N_2$ in panel
% (b) and the skewness $\gamma$ in panel (c).
Each panel shows the quartiles of the distribution of the estimate for
$100$ realisatons of time series generated by
\eqref{eq:snfnoise}. Each realisation was run until it reached $2000$
data points or $x=-1$ (indicating escape). Panel (a) shows the
estimate $c_\mathrm{emp}$ obtained by linear-least-squares fitting of
the approximate stationary Fokker-Planck equation~\eqref{eq:fp1} to
the empirical stationary density $p_\mathrm{emp}$, panel (b) shows the
quadratic coefficient $N_2$ obtained from \eqref{eq:fp2}, and panel
(c) shows the skewness $\gamma$ of $p_\mathrm{emp}$. For comparison,
we generate also $100$ linear time series using
\begin{equation}
  \label{eq:ou}
  \d x=-2\sqrt{\mu}x\d t+\d W_t
\end{equation}
and then fit $c_\mathrm{emp}$, $N_2$ and $\gamma$ for these time
series, too. Figure~\ref{fig:snftest} compares the quartiles of the
distributions for time series generated by the saddle-node normal form
\eqref{eq:snfnoise} and the linear equation \eqref{eq:ou}. If there is
small overlap in the distribution then the quantity is a good
indicator for the presence of a quadratic nonlinear term. We observe
that this is the case for all quanitites as long as the time series
has a length of $2000$ (see Figure~\ref{fig:snftest}(d) for the
distribution of time series lengths). Naturally, the uncertainty, and,
hence, the overlap, increases when the time series is shorter due to
escape from the potential well to $-\infty$ . This is the case for
smaller $\mu$ in Figure~\ref{fig:snftest} (see \citet{Thompson2011}
for a quantitative study of escape probabilities for small $\mu$). We
also notice that all quantities have a systematic deviation from the
theoretical value of the quantity they supposedly estimate: $N_2$
should be equal to $-1$ according to \eqref{eq:snfnoise}, the real
escape rate $c$ has a much smaller modulus than $c_\mathrm{emp}$
(compare Figure~\ref{fig:skewness}(a)), and the skewness $\gamma$ has
theoretically a much larger modulus (compare
Figure~\ref{fig:skewness}(b)). The deviation for $N_2$ is relatively
small, and is likely caused by the kernel density estimate for
$p_\mathrm{emp}$. The large modulus of $c_\mathrm{emp}$ shows that the
linear term $-\kappa_U\tilde x$ is a poor fit for $-\partial_xU(\tilde
x)$, which makes $c_\mathrm{emp}$ a measure of how much the odd
symmetry is broken by $-\partial_xU(\tilde x)$, rather than an
estimate for the escape rate. We note that $c_{\mathrm{emp},2}$ as
fitted in \eqref{eq:fp2} has a realistic modulus (close to zero, not
shown in Figure~\ref{fig:snftest}). The bias of the skewness $\gamma$
is mainly due to our restriction to time series that stay inside the
potential well. However, the characteristic dip of $\gamma$ seen in
Figure~\ref{fig:skewness}(b) is still visible in
Figure~\ref{fig:snftest}(c). \jsmod{\ref{sec:skewness}
  addresses the dependence of the empirical skewness on all method
  parameters.

  We also note that the estimates for the linear decay rate,
  $\kappa_U$ and $\kappa_\mathrm{ACF}$ are more accurate than the the
  estimates for the nonlinearity but by their nature they cannot
  distinguish the time series generated by the saddle-node normal form
  from a linear surrogate.}

The summative conclusion from Figure~\ref{fig:snftest} is that one can
detect the underlying nonlinearity of the deterministic part of
Equation~\eqref{eq:snfnoise} by observing $c_\mathrm{emp}$, $N_2$ or
the skewness $\gamma$ if the time series is moderately long.  One can
expect this to be true in the more general case of a time series
generated by a system with  deterministic dynamics close to a fold
and additive noise. While one can in principle recover the underlying
parameter $\mu$ from the estimates for $\kappa_U$, $N_2$ and $\sigma$
(see \citet{Thompson2011}), the uncertainty in $N_2$
propagates dramatically into uncertainty for $\mu$ (one has to divide
by $N_2$).

\jsmod{The results of Figure~\ref{fig:snftest} show the stationary
  case ($\epsilon=0$). For slowly drifting time series the same
  analysis will then have to be applied to parts of the time series in
  sliding windows. We demonstrate this for two geological times series
  in Section~\ref{sec:geolog}. For time series with rapidly drifting
  system parameter $\epsilon$ the estimates for the softening,
  $c_\mathrm{emp}$, $N_2$ and the skewness $\gamma$, will not give a
  detectable difference from the corresponding linear time series
  before tipping. The conclusion one would draw from this absence of
  detectable softening would be correct: tipping happens when the
  linear decay rate $\kappa$ reaches the critical value $0$ (or
  slightly later, see \citet{Thompson2011}).  The same holds if the
  noise level is low, which is equivalent to rapid drift (see
  \citet{Thompson2011} for the transformation). In practice one uses
  the argument the other way round: if the nonlinear estimates do not
  give a value significantly different from zero one will conclude that
  the ratio between noise level and drifting speed is small.}

\section{Studies of Geological Time Series}
\label{sec:geolog}

We now apply our nonlinear investigations to two ancient climate
tippings. The first is the major transition marking the end of the
latest ice age which occurred about $17,000$ years before the present
(yrs BP). The data used for this is a temperature reconstruction from
the Vostok ice core deuterium record (\citet{Petit1999}). The second
more recent tipping is the ending of the Younger Dryas using the
grayscale from the Cariaco basin sediments in Venezuela
(\citet{Hughen2000}). This Younger Dryas event was a curious cooling,
about $11,500$ years ago, just as the Earth was warming up after the
last ice age. It ended in a dramatic tipping point, when the Arctic
warmed by $7\celsius$ in $50$ years. This sudden ending has been
related (\citet{Houghton2004}) to a \emph{switching-on} of the global
oceanic thermohaline circulation (THC). This switch-on is known from
extensive theoretical studies (\citet{Dijkstra2005}) to be at a
saddle-node fold arising as a perturbation of a sub-critical pitchfork
(see, for example, \citet{Rahmstorf2000} and \citet{Thompson2011}).
\begin{figure}[ht]
  \centering
  \includegraphics[width=\textwidth]{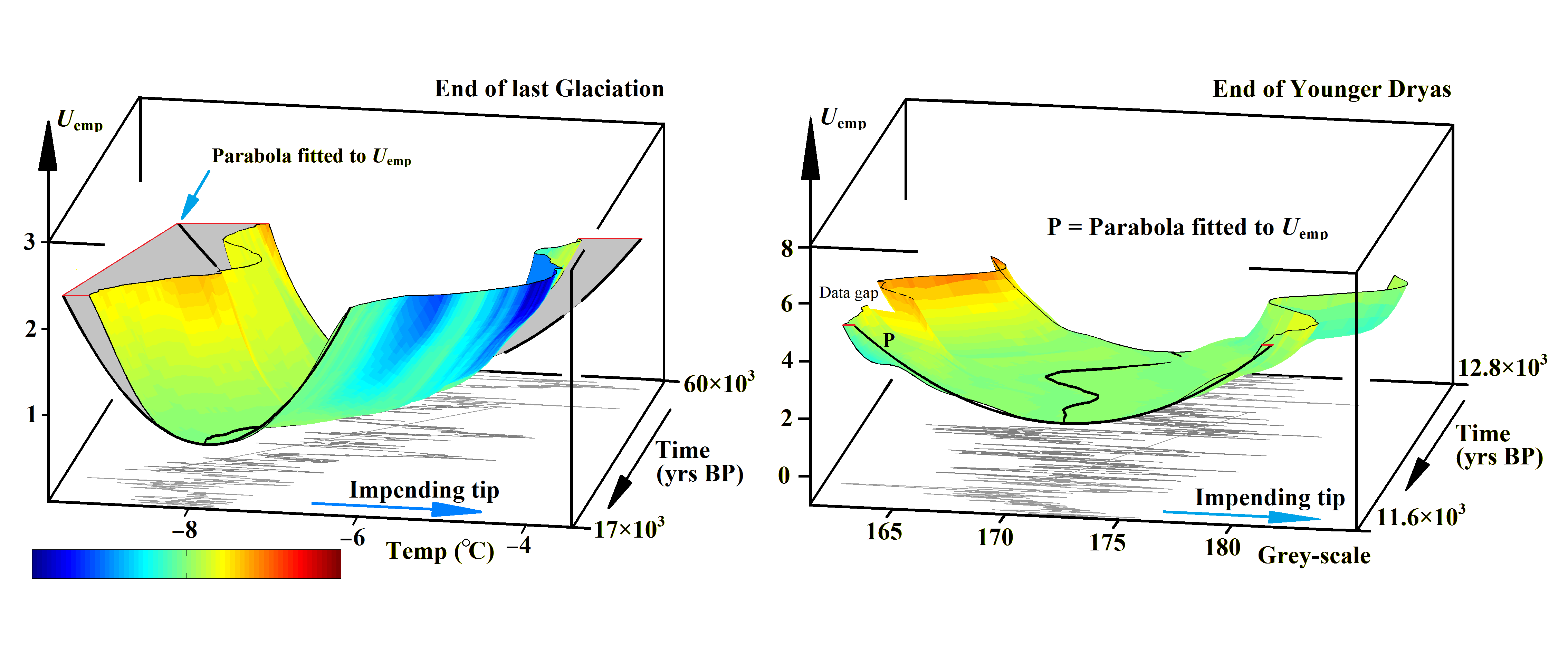}
  \caption{Two predicted potential energy surfaces for (a) the end of
    the last glaciation, using the Vostok ice-core record, and (b) the
    end of the Younger Dryas event when the Arctic warmed by
    $7\celsius$ in $50$ years, using the grey-scale of basin sediment
    in Cariaco, Venezuela. The colour shows the deviation from a
    (time-dependent) parabola fitted to $U_\mathrm{emp}(x,t)$, three
    of which are illustrated above. Red signifies a positive
    deviation, blue a negative deviation. Time is given in years
    before the present (BP).  \protect\jsmod{The plotted potential was
      obtained via
      $U_\mathrm{emp}=-\sigma_\mathrm{emp}^2\log(p_\mathrm{emp})/2$
      and \eqref{eq:sigmaemp} and sliding windows of half the length
      of the record.}}
  \label{fig:c366}
\end{figure}

\begin{figure}[th]
  \centering
  \includegraphics[width=0.8\textwidth]{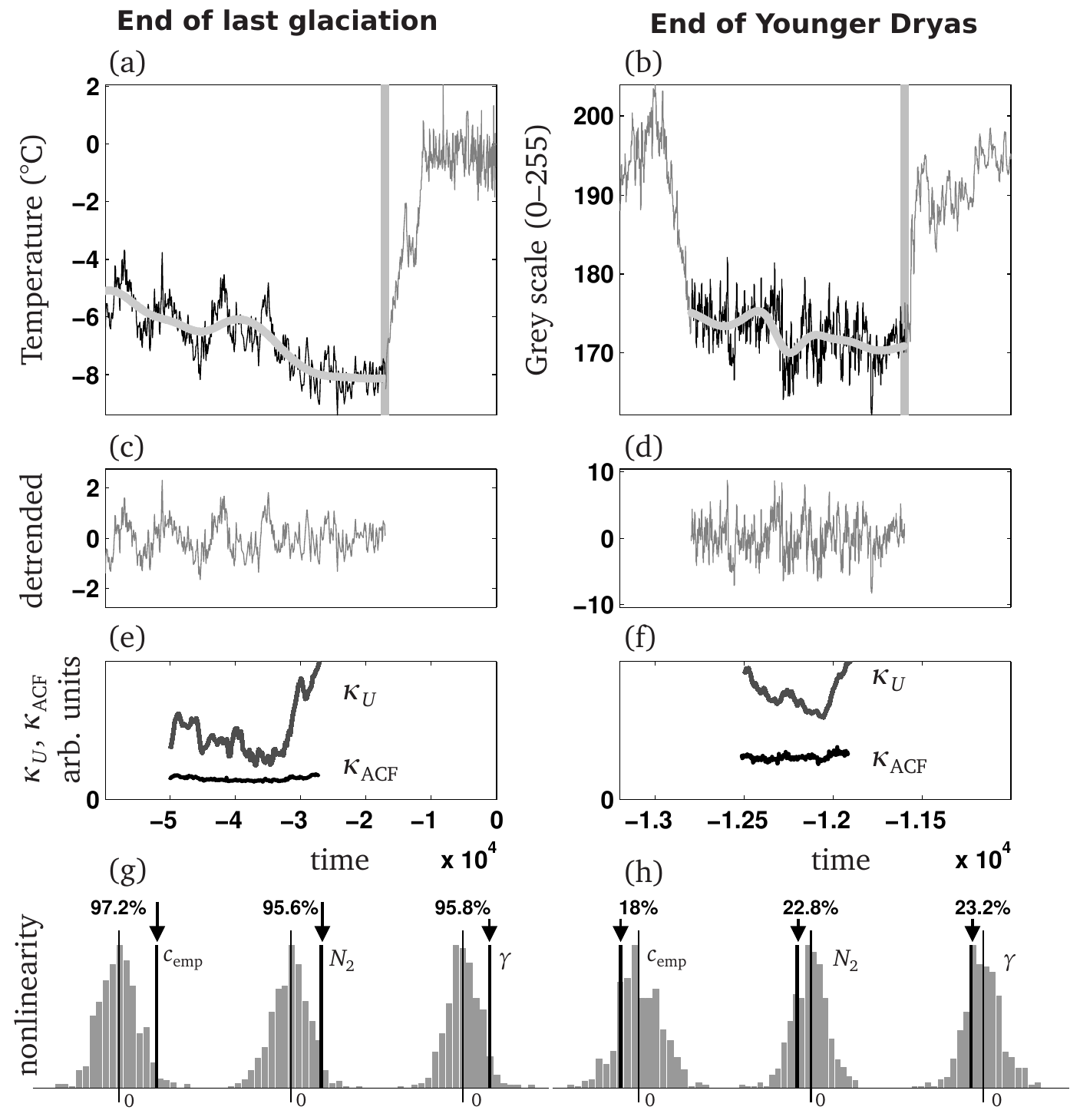}  
  \caption{Linear and nonlinear coefficients for two paleo-climate
    time series. (a,c,e,g): End of last glaciation (snapshot of data
    from \citet{Petit1999}); (b,d,f,h): End of Younger Dryas (snapshot
    of data from \citet{Hughen2000}). Only the black part of the data
    in panels (a) and (b) was used in the analysis. Panel (c) and (d)
    show the time series after detrending. Panel (e) and (f) show the
    linear indicators $\kappa_\mathrm{ACF}$ and $\kappa_U$. Panel (g)
    and (h) show the means of the estimates $c_\mathrm{emp}$, $N_2$
    and the skewness $\gamma$, compared to zero (thin vertical black
    line) and a histogram of estimates sampled from $500$ random
    linear time series generated by the linear process \eqref{eq:ou}
    with
    $-2\sqrt{\smash[b]{\mu}}=\kappa_\mathrm{ACF}$. \protect\jsmod{Sliding
      window length $w=N/2$.}}
  \label{fig:records}
\end{figure}

Under linear time-series analysis of the directly preceding data,
neither of these two events shows a strong trend in the stability
propagator (see Figure~\ref{fig:records}). Meanwhile sample estimates
of the underlying potential functions are shown in
Figure~\ref{fig:c366}.

In Figure~\ref{fig:c366}(a), for the end of the last glaciation, we
see clearly that the well is softening (falling beneath the parabolic
fit) on the high-temperature end while hardening (rising above the
parabola) on the low-temperature end. So there is a strong nonlinear
signal that a jump to higher temperatures may be pending (as indeed it
was). We show that this signal is statistically significant in our
full analysis in Figure~\ref{fig:records}. By comparison, the study of
the Younger Dryas, illustrated in Figure~\ref{fig:c366}(b), provided
no significant nonlinear conclusions.

% Figure~\ref{fig:records} shows two paleo-climate time series from
% periods prior to major climate
% transitions. Figure~\ref{fig:records}(a) shows a temperature
% reconstruction from just before the end of the last glaciation from
% the Vostok ice core deuterium record
% \citep{Petit1999}. Figure~\ref{fig:records}(b) shows greyscale records
% from the Cariaco basin \citep{Hughen2000} for a time span before the
% end of the Younger Dryas. We
For Figure~\ref{fig:records} we detrended both time series using
Gaussian filtering (see Figures~\ref{fig:records}(c) and (d)) and
estimated the linear indicators $\kappa_\mathrm{ACF}$ and
$\kappa_U$. We observe weak to non-existent trends, leaving the
evidence inconclusive at the linear level (grey and black curves
Figure~\ref{fig:records}(e) and (f)). The ratio between
$\kappa_\mathrm{ACF}$ and $\kappa_U$, both shown in
Figure~\ref{fig:records}(e) and (f), gives an estimate of the variance
$\sigma^2$ of the noise input. \jsmod{The unit for
  $\kappa_\mathrm{ACF}$ is \emph{units of $x$ per time step} (so, for
  example for Fig.~\ref{fig:records} (e) $\celsius/\Delta t$), and,
  correspondingly, the unit for $\kappa_U$ is $\Delta t/\celsius$
  (unit of $\kappa_\mathrm{ACF}$ divided by unit for $\sigma^2$). Both
  units are arbitrary such that we do not indicate the scale of the
  $y$-axis in Figure~\ref{fig:records}(e,f).}

\jsmod{\lenton discuss the linear analysis of both time
  series in greater detail. Their time series \emph{Vostok}
  corresponds to Figure~\ref{fig:records}(a), their time series
  \emph{Cariaco} corresponds to Figure~\ref{fig:records}(b). (Note
  that \lenton use the deuterium proxy directly whereas
  Figure~\ref{fig:records}(a) shows the temperature reconstruction.)
  Specifically, Figures 2 and 4, and Figures 9 and A1 in
  \lenton discuss how the results depend on method
  parameters such as detrending bandwidth and sliding window. Their
  analysis finds that the presence of a significant trend depends
  strongly on both parameters.  See \lenton, Table~1 (rows
  \emph{Vostok} and \emph{Cariaco}), for a summary of the evidence at
  the linear level.}

At the nonlinear level, the time series for the end of the last
glaciation has a strong quadratic nonlinearity in $\partial_xU(\tilde
x)$ where the potential well $U$ is presumed to govern the
deterministic part of the dynamics (Figure~\ref{fig:records}(g)). The
indicators $c_\mathrm{emp}$, $N_2$ and the skewness $\gamma$ in
Figure~\ref{fig:records}(g) are all far removed from what can be
expected by chance in a linear time series. The histogram in the
background of Figure~\ref{fig:records}(g,h) has been sampled from
$500$ random linear time series generated by the linear process given
in Equation \eqref{eq:ou} with
$-2\sqrt{\smash[b]{\mu}}=\kappa_\mathrm{ACF}$ where
$\kappa_\mathrm{ACF}$ is taken from the estimate shown in
Figure~\ref{fig:records}(e). The percentages at the top of
Figure~\ref{fig:records}(g,h) express how far in the tail of the
histogram the \jsmod{mean of the} quantity extracted from the
geological time series is. \jsrem{ $50$\% corresponds to the median of
  the histogram, a percentage smaller than $50$\% gives the percentage
  of linear realisations that were further away from the median than
  the geologocal data.} \jsmod{We warn that the three quantities
  $\gamma$, $N_2$ and $c_\mathrm{emp}$ are not really three
  independent indicators as they all depend on the same estimate of
  the empirical density $p_\mathrm{emp}$.}

Visual inspection of the well shape in Figure~\ref{fig:c366} confirms
that the well is softening (bending downward) on the high-temperature
end and hardening (bending upward) on the low-temperature end. So, at
the nonlinear level the time series data close to the bottom of the
well gives already evidence for a propensity to escape toward larger
temperatures. The time series for the end of the Younger Dryas does
not show evidence for strong nonlinearity of the underlying dynamics
at the second-order level (Figure~\ref{fig:records}(h)). The values
for $c_\mathrm{emp}$, $N_2$ and the skewness $\gamma$ can all be
explained by randomness as the histograms of estimates for the $500$
linear time series (as the histograms in Figure~\ref{fig:records}(h)
show). Note that we did not include the previous transition (visible
at the left end in Figure~\ref{fig:records}(b)) into our analysis as
our aim was to infer the nonlinearity exclusively from the data near
the tipping equilibrium.

\section{Conclusion}
\label{sec:conc}
\jsmod{The main message of the paper is that the linear analysis
  investigated by \lenton on its own, even if all estimates
  are accurate, does not contain all information necessary to estimate
  the probability of tipping over time.} The result of the linear
analysis is an estimate of the decay rate $\kappa$\jsmod{, which is taken}
relative to the time spacing of the measurements. Even if this decay
rate has an identifiable trend to zero the probability of tipping over
time is determined by the dominant nonlinear coefficient $N_2$ in
conjunction with $\kappa$ and the noise level.

We found that the dominant nonlinear coefficient $N_2$ is much harder
to estimate accurately, so we also looked at proxies that indicate
nonlinearity such as the skewness $\gamma$ (as proposed by
\citet{Guttal2008}) or the coefficient $c_\mathrm{emp}$ from
Equation~\eqref{eq:fp1} (which is also a scalar signed measure for
deviation from linear theory). Our study of the saddle-node normal
form suggests that both proxies, $\gamma$ and $c_\mathrm{emp}$, give
values that are easier to distinguish from random chance than the
estimate for $N_2$ itself. A significantly non-zero value of either of
these proxies indicates that a quadratic nonlinear term is present and
gives its sign. 

However, we found the uncertainty for moderately long time series
($N=2000$) too large to translate the proxy back into a quantitatively
reliable estimate for the normal form parameter (this would be
necessary to read off the probability for tipping from the tables in
\citet{Thompson2011}). We also found that the nonlinear proxies do not
have discernible trends when one approaches tipping, so it makes sense
only to extract their overall mean from the time series but not the
trend.
\bibliographystyle{unsrtnat}
%\bibliography{climate}

\begin{figure}[ht]
  \centering
  \includegraphics[width=\textwidth]{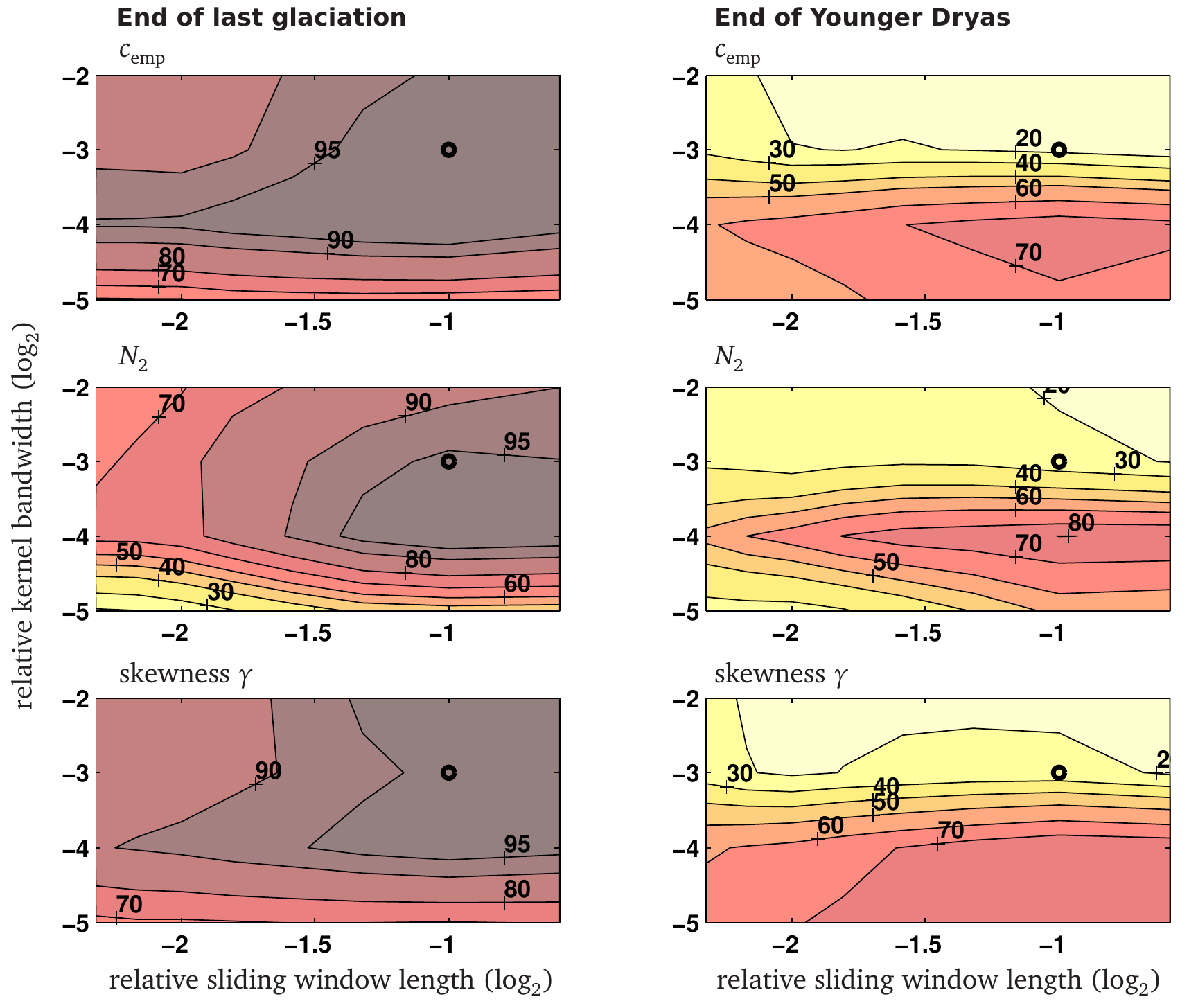}
  \caption{Sensitivity of estimates with respect to the method
    parameters sliding window length (horizontal) and fitting kernel
    bandwidth (vertical). The contour levels show every $10$th
    percentile ($10\%,\ldots,90\%$) and the $95\%$ level (colour online).}
  \label{fig:sens}
\end{figure}

%\cbstart
\appendix{Sensitivity analysis for estimates of nonlinear quantities
  in palaeoclimate records}
\label{sec:sensitivity}
Figure~\ref{fig:sens} shows how the results in
Figure~\ref{fig:records}(g,h) depend on our method parameters.
Specifically, we vary the length of the sliding window and the
bandwidth of the Gaussian kernel when fitting the quasi-equilibrium
drift.
In all panels the $x$-axis is the length of the sliding window
relative to the overall length of the time series on a logarithmic
scale (overall length is $N=525$ for the end of the last glaciation,
$N=2646$ for the end of the Younger Dryas). For example, the $x$-value
of $-1$ corresponds to a window length $w=N/2$. The $y$-axis is the
kernel bandwidth for the Gaussian fitting kernel relative to the
overall length of the time series on a logarithmic scale (this enters
the function \texttt{kde1d} as a resolution parameter). For example,
the $y$-value of $-3$ means that the resolution parameter of
\texttt{kde1d} is $2^3=8$.

For each measured quantity we compare the mean of its value over all
sliding windows to the percentiles of the linear surrogate. In this
way contour levels close to $0$ or $100$ correspond to significant
deviations from a time series generated by a linear process. The
distributions of these time series are shown as gray background in
Fig.~\ref{fig:records}(g,h). The value shown as a black vertical line
in Figure~\ref{fig:records}(g,h) was obtained using the method
parameters highlighted by a black circle in Figure~\ref{fig:sens}.

Fig.~\ref{fig:sens} gives evidence that the nonlinear features of the
time series for the end of the last glaciation are robust as long as
the detrending remains coarse grained. The end of the Younger Dryas
does not show any nonlinear feature that is significantly different
from the surrogates for any choise of method parameters.
\appendix{Dependence of skewness on distance to tipping point}
\label{sec:skewness}
The observations in Figure~\ref{fig:skewness}(b) and
Figure~\ref{fig:snftest}(c) appear to contradict the statement by
\citet{Guttal2008} that skewness increases as one approaches tipping.
Furthermore, there is a quantiative difference between
Figure~\ref{fig:skewness}(b) and Figure~\ref{fig:snftest}(c), even
though they show the same qualitative feature: non-monotonicity of the
skewness in the control parameter $\mu$.  Note that \citet{Guttal2008}
studied the same type of tipping as Figure~\ref{fig:skewness},
however, not in the fold normal form but in a particular ecological
model. In order to investigate this apparent contradiction we look at
how the skewness depends on two parameters: the control parameter
$\mu$ and the value at which we consider a realisation as
escaped. More precisely, both figures, \ref{fig:skewness}(b) and
\ref{fig:snftest}(c), study
\begin{equation}
  \label{eq:snfapp}
  \d x=\left[\mu-x^2\right]\d t+\sigma \d W_t\mbox{,}
\end{equation}
where $\sigma=1$ without loss of generality and $\mu$ is constant. Any
realisation of \eqref{eq:snfapp} escapes to $-\infty$ in finite time
almost surely such that \eqref{eq:snfapp} does not possess a
stationary density. In the calculation for Figure~\ref{fig:skewness}
we modified the process \eqref{eq:snfapp} by re-injecting the
realisation at $+\infty$ whenever it escaped to $-\infty$. In this way
the process has a stationary density given by the Fokker-Planck
equation \eqref{eq:fp} with $-\partial_xU(x)=\mu-x^2$. However, this
is not accurately describing the density that one is interested in:
the probability density of the realisations under the condition that
the realisation has not yet `escaped'. Clearly, this conditional
probability depends on what one means by escaped. So, we introduce the
escape boundary $b$ as an additional parameter. Whenever, a
realisation gets smaller than the value $-\sqrt{\mu}-b$ (that is, it
has run away beyond the local maximum of the potential by distance
$b$) we stop the integration of \eqref{eq:snfapp}. When one considers
the set of all realisations before escape then this also leads to a
stationary density. This was done for Figure~\ref{fig:snftest} using
$b=-1+\sqrt{\mu}$ (such that we count a realisation as escaped
whenever it reaches $-1$). The re-injection used for
Figure~\ref{fig:skewness} and the choice of $x=-1$ as escape indicator
lead to the systematic difference in skewness $\gamma$ between the
figures \ref{fig:snftest}(c) and \ref{fig:skewness}(b). In practice
one encounters only the conditional probability as shown in
Figure~\ref{fig:snftest}. So, it is important to check how the
skewness depends on the somewhat arbitrary choice of the the value $b$
indicating escape. Moreover, the estimate shown in
Figure~\ref{fig:snftest}(c) suffers from the shorteness of the time
series for $\mu$ close to $0$ (which is typical in
practice). 

\begin{figure}[ht]
  \centering
  \includegraphics[width=0.6\textwidth]{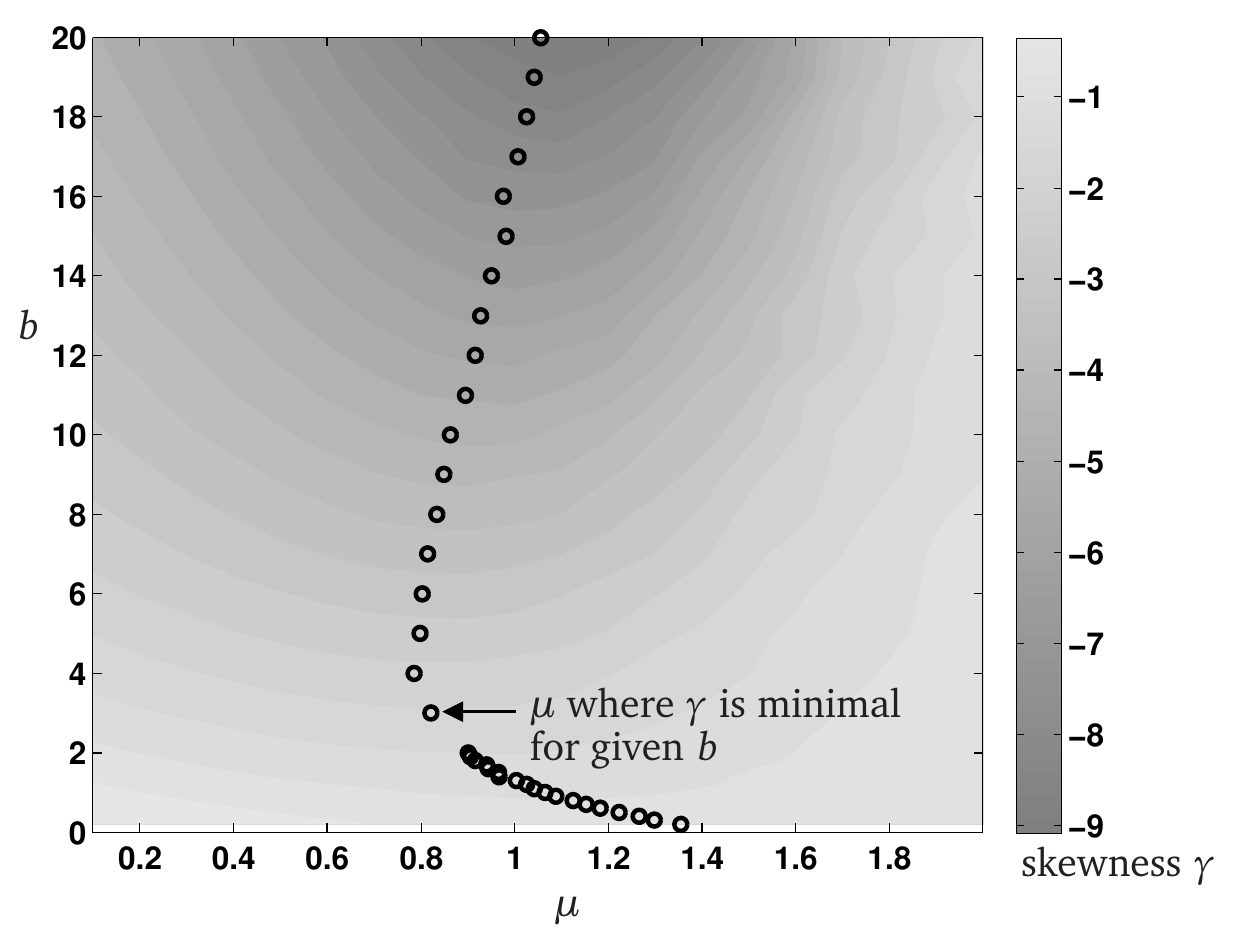}
  \caption{Skewness of probability density of non-escaped realisations}
  \label{fig:skewscan}
\end{figure}
Figure~\ref{fig:skewscan} shows the result of a more precise
calculation: we simulate \eqref{eq:snfapp} for a large ensemble of
realisations ($n=100\,000$). At every time step $t_k$ we collect the
$n_{\mathrm{esc},k}$ realisations that have reached a value
$x\leq-\sqrt{\mu}-b$. They count as escaped and have to be
discarded. In order to avoid depletion of the ensemble we start
$n_{\mathrm{esc},k}$ new realisations at $t_k$. For each newly started
realisation we pick the current value of one randomly chosen
non-escaped realisation as the inital value. This process reaches a
stationary density, which is an approximation of the conditional
probability density of $x$ under the condition that $x$ has not yet
escaped beyond $-\sqrt{\mu}-b$. We observe in
Figure~\ref{fig:skewscan} that the depth of the dip in the skewness
depends on $b$ (indicated by grayscale) but that the location of the
dip is uniformly larger than the critical value $0$ of the control
parameter $\mu$.

It is possible that the simulations by \citet{Guttal2008} stayed
consistently to the right of the curve of extreme skewness shown in
Figure~\ref{fig:skewscan} such that only the increase is observed. Note
that in their model the notion of skewness has the opposite sign, so
one should consider the modulus of the skewness.
%\cbend
\end{document}